\documentclass[a4paper, 12 pt]{article}

\usepackage[T2A]{fontenc}
\usepackage[cp1251]{inputenc}
\usepackage[english]{babel}
\usepackage[tbtags]{amsmath}
\usepackage{amsfonts,amssymb}

\usepackage{mathrsfs}

\usepackage{geometry} 
\begin{document}
\title{Еxponential Riordan arrays and generalized Narayana polynomials}
 \author{E. Burlachenko}
 \date{}

 \maketitle
\begin{abstract}
Generalized Euler polynomials ${{\alpha }_{n}}\left( x \right)={{\left( 1-x \right)}^{n+1}}\sum\nolimits_{m=0}^{\infty }{{{p}_{n}}}\left( m \right){{x}^{m}}$, where ${{p}_{n}}\left( x \right)$ is the polynomial of degree $n$, are the numerator polynomials  of the generating functions of diagonals of the ordinary Riordan arrays. Generalized Narayana polynomials ${{\varphi }_{n}}\left( x \right)={{\left( 1-x \right)}^{2n+1}}\sum\nolimits_{m=0}^{\infty }{\left( m+1 \right)...\left( m+n \right){{p}_{n}}}\left( m \right){{x}^{m}}$ are the numerator polynomials  of the generating functions of diagonals of the exponential Riordan arrays. In present paper we consider the constructive relationship between these two types of numerator polynomials.
\end{abstract}
\section{Introduction}
This paper is a continuation of the paper “Riordan arrays and generalized Euler polynomials” [1]. In Section 2 we will  briefly retell its content. For the integrity of presentation, we will change some notation adopted in [1]. 

Subject of our study is the transformations in space of formal power series and the corresponding matrices. We associate rows and columns of matrices with the generating functions of their elements. $n$th coefficient of the series $a\left( x \right)$, $n$th row, $n$th descending diagonal and $n$th column of the matrix $A$ will be denoted  respectively by
$$\left[ {{x}^{n}} \right]a\left( x \right),\qquad \left[ n,\to  \right]A,   \qquad\left[ n,\searrow\right ]A,   \qquad A{{x}^{n}}.$$
Matrix $\left( f\left( x \right),g\left( x \right) \right)$, ${{g}_{0}}=0$, $n$th column of which, $n=0,\text{ }1,\text{ }2,\text{ }...$ , has the generating function $f\left( x \right){{g}^{n}}\left( x \right)$, is called Riordan array [2] – [6]. It is the product of two matrices that correspond to multiplication and composition of series:
$$\left( f\left( x \right),g\left( x \right) \right)=\left( f\left( x \right),x \right)\left( 1,g\left( x \right) \right),$$
$$\left( f\left( x \right),x \right)a\left( x \right)=f\left( x \right)a\left( x \right), \qquad\left( 1,g\left( x \right) \right)a\left( x \right)=a\left( g\left( x \right) \right),$$
$$\left( f\left( x \right),g\left( x \right) \right)\left( b\left( x \right),a\left( x \right) \right)=\left( f\left( x \right)b\left( g\left( x \right) \right),a\left( g\left( x \right) \right) \right).$$
Matrices 
$${{\left| {{e}^{x}} \right|}^{-1}}\left( f\left( x \right),g\left( x \right) \right)\left| {{e}^{x}} \right|={{\left( f\left( x \right),g\left( x \right) \right)}_{{{e}^{x}}}},$$
where $\left| {{e}^{x}} \right|$ is the diagonal matrix, $\left| {{e}^{x}} \right|{{x}^{n}}={{{x}^{n}}}/{n!}\;$, are called exponential Riordan arrays. Denote
$$\left[ n,\to  \right]{{\left( f\left( x \right),g\left( x \right) \right)}_{{{e}^{x}}}}={{s}_{n}}\left( x \right), \qquad{{f}_{0}}\ne 0, \qquad{{g}_{1}}\ne 0.$$
Then
$${{\left( f\left( x \right),g\left( x \right) \right)}_{{{e}^{x}}}}{{\left( 1-\varphi x \right)}^{-1}}={{\left| {{e}^{x}} \right|}^{-1}}\left( f\left( x \right),g\left( x \right) \right){{e}^{\varphi x}}={{\left| {{e}^{x}} \right|}^{-1}}f\left( x \right)\exp \left( \varphi g\left( x \right) \right),$$
or
$$\sum\limits_{n=0}^{\infty }{\frac{{{s}_{n}}\left( \varphi  \right)}{n!}{{x}^{n}}}=f\left( x \right)\exp \left( \varphi g\left( x \right) \right).$$
Sequence of polynomials ${{s}_{n}}\left( x \right)$ is called Sheffer sequence, and in the case $f\left( x \right)=1$, binomial sequence. The properties of the Sheffer sequences are  the subject of study of the umbral calculus [7]. Matrix
$$P=\left( \frac{1}{1-x},\frac{x}{1-x} \right)={{\left( {{e}^{x}},x \right)}_{{{e}^{x}}}}=\left( \begin{matrix}
   1 & 0 & 0 & 0 & \cdots   \\
   1 & 1 & 0 & 0 & \cdots   \\
   1 & 2 & 1 & 0 & \cdots   \\
   1 & 3 & 3 & 1 & \cdots   \\
   \vdots  & \vdots  & \vdots  & \vdots  & \ddots   \\
\end{matrix} \right)$$
is called Pascal matrix. Power of the Pascal matrix is defined by the identity
$${{P}^{\varphi }}=\left( \frac{1}{1-\varphi x},\frac{x}{1-\varphi x} \right)={{\left( {{e}^{\varphi x}},x \right)}_{{{e}^{x}}}}.$$

Along with the lower triangular  Riordan  matrices, we will consider “square” matrices $\left( b\left( x \right),a\left( x \right) \right)$,  ${{b}_{0}}\ne 0$, ${{a}_{0}}=1$. For example,
$$\left( 1,\frac{1}{1+x} \right)=\left( \begin{matrix}
   1 & 1 & 1 & 1 & \cdots   \\
   0 & -1 & -2 & -3 & \cdots   \\
   0 & 1 & 3 & 6 & \cdots   \\
   0 & -1 & -4 & -10 & \cdots   \\
   \vdots  & \vdots  & \vdots  & \vdots  & \ddots   \\
\end{matrix} \right)$$
This includes the upper triangular matrix $\left( 1,1+x \right)$, whose transpose is the Pascal matrix and which coincides with the  matrix of shift operator:
$$\left( 1,1+x \right)={{P}^{T}}=E=\left( \begin{matrix}
   1 & 1 & 1 & 1 & \cdots   \\
   0 & 1 & 2 & 3 & \cdots   \\
   0 & 0 & 1 & 3 & \cdots   \\
   0 & 0 & 0 & 1 & \cdots   \\
   \vdots  & \vdots  & \vdots  & \vdots  & \ddots   \\
\end{matrix} \right).$$
Matrix $\left( b\left( x \right),a\left( x \right) \right)$ can be multiplied from the right by the matrix with the finite columns and from the left by the matrix with the finite rows. At first (before Section 4), we restrict ourselves to the set of matrices of the form $\left( 1,a\left( x \right) \right)$. Since
$$\left[ n,\to  \right]\left( 1,a\left( x \right) \right)=[n,\searrow ]\left( 1,xa\left( x \right) \right),$$
then the matrix $\left( 1,a\left( x \right) \right)$ is a tool for study of the matrix $\left( 1,xa\left( x \right) \right)$. Denote
$$\left[ n,\to  \right]\left( 1,a\left( x \right)-1 \right)={{v}_{n}}\left( x \right)=\sum\limits_{m=1}^{n}{{{v}_{m}}{{x}^{m}}}, \qquad n>0.$$
Since
$$\left( 1,a\left( x \right)-1 \right)\left( 1,1+x \right)=\left( 1,a\left( x \right) \right),  \qquad\left[ n,\to  \right]\left( 1,1+x \right)=\frac{{{x}^{n}}}{{{\left( 1-x \right)}^{n+1}}},$$
then
$$\left[ n,\to  \right]\left( 1,a\left( x \right) \right)=\sum\limits_{m=1}^{n}{\frac{{{v}_{m}}{{x}^{m}}}{{{\left( 1-x \right)}^{m+1}}}}=\sum\limits_{m=1}^{n}{\frac{{{v}_{m}}{{x}^{m}}{{\left( 1-x \right)}^{n-m}}}{{{\left( 1-x \right)}^{n+1}}}=}\frac{{{\alpha }_{n}}\left( x \right)}{{{\left( 1-x \right)}^{n+1}}}.$$
If $a\left( x \right)={{e}^{x}}$, then ${{\alpha }_{n}}\left( x \right)={{{A}_{n}}\left( x \right)}/{n!}\;$, where ${{A}_{n}}\left( x \right)$ are the Euler polynomials:
$$\frac{{{A}_{n}}\left( x \right)}{{{\left( 1-x \right)}^{n+1}}}=\sum\limits_{m=0}^{\infty }{{{m}^{n}}}{{x}^{m}},  \qquad{{A}_{n}}\left( 1 \right)=n!.$$
For example,
$${{A}_{1}}\left( x \right)=x,  \qquad{{A}_{2}}\left( x \right)=x+{{x}^{2}},  \qquad{{A}_{3}}\left( x \right)=x+4{{x}^{2}}+{{x}^{3}},$$
$${{A}_{4}}\left( x \right)=x+11{{x}^{2}}+11{{x}^{3}}+{{x}^{4}}.$$
In this connection we will called these polynomials the generalized Euler polynomials (GEP).

“Square” Riordan arrays (called convolution arrays) and numerator polynomials of the generating functions of their rows were considered in the series of papers [8] – [12]. In [13] such matrices are called generalized Riordan arrays. Concept of generalized Euler polynomials (called ${{p}_{n}}$-associated Eulerian polynomials)  in general form is represented in [14]. 

Denote 
$$\left[ n,\searrow  \right]\left( 1,xa\left( x \right) \right)=\frac{{{\alpha }_{n}}\left( x \right)}{{{\left( 1-x \right)}^{n+1}}},  \qquad\left[ n,\to  \right]{{\left( 1,\log a\left( x \right) \right)}_{{{e}^{x}}}}={{u}_{n}}\left( x \right),$$
$$\left[ n,\to  \right]\left( 1,a\left( x \right)-1 \right)={{v}_{n}}\left( x \right).$$
If the sequence of polynomials has the form ${{c}_{0}}\left( x \right)=1$, $\left[ {{x}^{0}} \right]{{c}_{n}}\left( x \right)=0$, we will bear in mind that the expression $\left( {1}/{x}\; \right){{c}_{n}}\left( x \right)$ corresponds to the case $n>0$. Denote
$$\frac{1}{x}{{\alpha }_{n}}\left( x \right)={{\tilde{\alpha }}_{n}}\left( x \right),  \quad\frac{1}{x}{{A}_{n}}\left( x \right)={{\tilde{A}}_{n}}\left( x \right),  \quad\frac{1}{x}{{u}_{n}}\left( x \right)={{\tilde{u}}_{n}}\left( x \right), 
\quad\frac{1}{x}{{v}_{n}}\left( x \right)={{\tilde{v}}_{n}}\left( x \right).$$
Let the symbols ${{\left( \varphi  \right)}_{n}}$, ${{\left[ \varphi  \right]}_{n}}$ denote respectively the falling and the rising  factorial:
$${{\left( \varphi  \right)}_{n}}=\varphi \left( \varphi -1 \right)...\left( \varphi -n+1 \right),   \qquad{{\left[ \varphi  \right]}_{n}}=\varphi \left( \varphi +1 \right)...\left( \varphi +n-1 \right).$$
In Section 2 we consider the generalized Euler polynomials and associated transformations. We introduce the matrices ${{\tilde{U}}_{n}}$, ${{\tilde{V}}_{n}}$:
$${{\tilde{U}}_{n}}{{x}^{p}}=\frac{1}{n!}{{\left( 1-x \right)}^{n-1-p}}{{\tilde{A}}_{p+1}}\left( x \right),   \qquad\tilde{U}_{n}^{-1}{{x}^{p}}={{\left( x-1 \right)}_{p}}{{\left[ x+1 \right]}_{n-p-1}},$$
$${{\tilde{V}}_{n}}{{x}^{p}}={{\left( 1+x \right)}^{n-p-1}}{{x}^{p}},    \qquad\tilde{V}_{n}^{-1}{{x}^{p}}={{\left( 1-x \right)}^{n-p-1}}{{x}^{p}},  \qquad p=0, \text{ }1,\text{ } … , \text{ } n-1.$$
Then
$${{\tilde{U}}_{n}}{{\tilde{u}}_{n}}\left( x \right)={{\tilde{\alpha }}_{n}}\left( x \right),   \qquad{{\tilde{V}}_{n}}{{\tilde{\alpha }}_{n}}\left( x \right)={{\tilde{v}}_{n}}\left( x \right).$$
We consider the series $_{\left( \beta  \right)}a\left( x \right)$, $_{\left( 0 \right)}a\left( x \right)=a\left( x \right)$, that are defined as follows:
$$_{\left( \beta  \right)}{{a}^{\varphi }}\left( x \right)=\sum\limits_{n=0}^{\infty }{\frac{\varphi }{\varphi +n\beta }}\frac{{{u}_{n}}\left( \varphi +n\beta  \right)}{n!}{{x}^{n}}.$$
Denote
$$\left[ n,\searrow  \right]\left( 1,x{}_{\left( \beta  \right)}a\left( x \right) \right)=\frac{_{\left( \beta  \right)}{{\alpha }_{n}}\left( x \right)}{{{\left( 1-x \right)}^{n+1}}},  \qquad\frac{1}{x}{}_{\left( \beta  \right)}{{\alpha }_{n}}\left( x \right)={}_{\left( \beta  \right)}{{\tilde{\alpha }}_{n}}\left( x \right).$$
We introduce the matrices 
$$A_{n}^{\beta }={{\tilde{U}}_{n}}{{E}^{n\beta }}\tilde{U}_{n}^{-1}=\tilde{V}_{n}^{-1}\tilde{D}{{\left( {{\left( 1+x \right)}^{n\beta }},x \right)}^{T}}{{\tilde{D}}^{-1}}{{\tilde{V}}_{n}}, \qquad\tilde{D}{{x}^{n}}=\left( n+1 \right){{x}^{n}}.$$
Then
$$A_{n}^{\beta }{{\tilde{\alpha }}_{n}}\left( x \right)={}_{\left( \beta  \right)}{{\tilde{\alpha }}_{n}}\left( x \right).$$
 We give a general formula for the GEP associated with the generalized binomial series. Namely, let
$$_{\left( \beta  \right)}{{a}^{\varphi }}\left( x \right)=\sum\limits_{n=0}^{\infty }{\frac{\varphi }{\varphi +n\beta }}\left( \begin{matrix}
   \varphi +n\beta   \\
   n  \\
\end{matrix} \right){{x}^{n}},  \quad\frac{_{\left( \beta  \right)}{{\alpha }_{n}}\left( x \right)}{{{\left( 1-x \right)}^{n+1}}}=\sum\limits_{m=0}^{\infty }{\frac{m}{m+n\beta }}\left( \begin{matrix}
   m+n\beta   \\
   n  \\
\end{matrix} \right){{x}^{m}}.$$
Then
$$_{\left( \beta  \right)}{{\alpha }_{n}}\left( x \right)=\frac{1}{n}\sum\limits_{m=1}^{n}{\left( \begin{matrix}
   n\left( 1-\beta  \right)  \\
   m-1  \\
\end{matrix} \right)\left( \begin{matrix}
   n\beta   \\
   n-m  \\
\end{matrix} \right){{x}^{m}}}.$$

In Section 3 we consider the generalized Narayana polynomials ${{\varphi }_{n}}\left( x \right)$, which are the numerator polynomials of the matrix ${{\left( 1,xa\left( x \right) \right)}_{{{e}^{x}}}}$: 
$$\left[ n,\searrow  \right]{{\left( 1,xa\left( x \right) \right)}_{{{e}^{x}}}}=\frac{{{\varphi }_{n}}\left( x \right)}{{{\left( 1-x \right)}^{2n+1}}}=\sum\limits_{m=0}^{\infty }{\frac{{{\left[ m+1 \right]}_{n}}{{u}_{n}}\left( m \right)}{n!}{{x}^{m}}}, \quad\frac{1}{x}{{\varphi }_{n}}\left( x \right)={{\tilde{\varphi }}_{n}}\left( x \right).$$
We introduce the matrices ${{\tilde{F}}_{n}}$:
$${{\tilde{F}}_{n}}{{x}^{p}}={{\left( 1-x \right)}^{2n+1}}\sum\limits_{m=1}^{\infty }{{{m}^{p+1}}}\left( \begin{matrix}
   m+n  \\
   n  \\
\end{matrix} \right){{x}^{m-1}},$$
$$\tilde{F}_{n}^{-1}{{x}^{p}}=\frac{n!}{\left( 2n \right)!}{{\left( x-1 \right)}_{p}}{{\left[ x+n+1 \right]}_{n-p-1}},  \qquad p=0, \text{ }1,\text{ } … , \text{ } n-1.$$
Then
$${{\tilde{F}}_{n}}{{\tilde{u}}_{n}}\left( x \right)={{\tilde{\varphi }}_{n}}\left( x \right).$$
We introduce the matrices ${{\tilde{S}}_{n}}={{\tilde{F}}_{n}}\tilde{U}_{n}^{-1}$. Then
$${{\tilde{S}}_{n}}{{\tilde{\alpha }}_{n}}\left( x \right)={{\tilde{\varphi }}_{n}}\left( x \right).$$
It turns out that
$${{\tilde{S}}_{n}}=\tilde{V}_{n}^{-1}{{\tilde{C}}_{n}}{{\tilde{V}}_{n}},  \qquad{{\tilde{C}}_{n}}{{x}^{p}}=\frac{\left( n+p+1 \right)!}{\left( p+1 \right)!}{{x}^{p}}.$$
 We give a general formula for the GNP associated with the generalized binomial series. Namely, let
$$\frac{_{\left( \beta  \right)}{{\varphi }_{n}}\left( x \right)}{{{\left( 1-x \right)}^{2n+1}}}=\sum\limits_{m=0}^{\infty }{\frac{m}{m+\beta n}}\left( \begin{matrix}
   m+\beta n  \\
   n  \\
\end{matrix} \right){{\left[ m+1 \right]}_{n}}{{x}^{m}}.$$
Then
$$_{\left( \beta  \right)}{{\varphi }_{n}}\left( x \right)=\frac{\left( n+1 \right)!}{n}\sum\limits_{m=1}^{n}{\left( \begin{matrix}
   n\left( 2-\beta  \right)  \\
   m-1  \\
\end{matrix} \right)\left( \begin{matrix}
   n\beta   \\
   n-m  \\
\end{matrix} \right){{x}^{m}}}.$$

In Section 4 we consider transformations of the general form. Let ${{g}_{n}}\left( x \right)$, ${{h}_{n}}\left( x \right)$ are the numerator polynomialsof of the matrices $\left( b\left( x \right),xa\left( x \right) \right)$, ${{\left( b\left( x \right),xa\left( x \right) \right)}_{{{e}^{x}}}}$, ${{b}_{0}}\ne 0$, respectively. Denote
$$\left[ n,\to  \right]{{\left( b\left( x \right),\log a\left( x \right) \right)}_{{{e}^{x}}}}={{s}_{n}}\left( x \right).$$
We introduce the matrices  ${{U}_{n}}$, ${{F}_{n}}$:
$${{U}_{n}}{{x}^{p}}={{\left( 1-x \right)}^{n+1}}\frac{1}{n!}\sum\limits_{m=0}^{\infty }{{{m}^{p}}}{{x}^{m}},  \qquad U_{n}^{-1}{{x}^{p}}={{\left( x \right)}_{p}}{{\left[ x+1 \right]}_{n-p}};$$
$${{F}_{n}}{{x}^{p}}={{\left( 1-x \right)}^{2n+1}}\sum\limits_{m=0}^{\infty }{{{m}^{p}}}\left( \begin{matrix}
   m+n  \\
   n  \\
\end{matrix} \right){{x}^{m}},  \quad F_{n}^{-1}{{x}^{p}}=\frac{n!}{\left( 2n \right)!}{{\left( x \right)}_{p}}{{\left[ x+n+1 \right]}_{n-p}},$$
$p=0$,  $1$, … , $n$. Then
$${{U}_{n}}{{s}_{n}}\left( x \right)={{g}_{n}}\left( x \right),  \qquad{{F}_{n}}{{s}_{n}}\left( x \right)={{h}_{n}}\left( x \right).$$
We introduce the matrices ${{S}_{n}}={{F}_{n}}U_{n}^{-1}$. Then
$${{S}_{n}}{{g}_{n}}\left( x \right)={{h}_{n}}\left( x \right).$$
It turns out that
$${{S}_{n}}=V_{n}^{-1}{{C}_{n}}{{V}_{n}},$$  
$${{V}_{n}}{{x}^{p}}={{\left( 1+x \right)}^{n-p}}{{x}^{p}},    \quad V_{n}^{-1}{{x}^{p}}={{\left( 1-x \right)}^{n-p}}{{x}^{p}},  \quad{{C}_{n}}{{x}^{p}}=\frac{\left( n+p \right)!}{p!}{{x}^{p}} ;$$
$${{S}_{n}}{{x}^{p}}=\frac{\left( n+p \right)!\left( n-p \right)!}{n!}\sum\limits_{m=p}^{n}{\left( \begin{matrix}
   n  \\
   m-p  \\
\end{matrix} \right)\left( \begin{matrix}
   n  \\
   n-m  \\
\end{matrix} \right){{x}^{m}}},$$
$$S_{n}^{-1}{{x}^{p}}=\frac{p!\left( n-p \right)!}{\left( 2n \right)!}\sum\limits_{m=p}^{n}{\left( \begin{matrix}
   -n  \\
   m-p  \\
\end{matrix} \right)\left( \begin{matrix}
   2n  \\
   n-m  \\
\end{matrix} \right){{x}^{m}}}.$$

In Section 5 we consider the generalized Narayana polynomials of type B, which are the numerator polynomials of the matrix ${{\left( a\left( x \right),xa\left( x \right) \right)}_{{{e}^{x}}}}$, and similar polynomials, which are the numerator polynomials of the matrix $\left( {{\left( xa\left( x \right) \right)}^{\prime }},xa\left( x \right) \right)$.

In Section 6 we return to the series $_{\left( \beta  \right)}a\left( x \right)$ from Section 2 and consider the transformations
$$G_{n}^{\beta }={{U}_{n}}{{E}^{n\beta }}U_{n}^{-1}=V_{n}^{-1}{{\left( {{\left( 1+x \right)}^{n\beta }},x \right)}^{T}}{{V}_{n}},$$
$$H_{n}^{\beta }={{F}_{n}}{{E}^{n\beta }}F_{n}^{-1}={{S}_{n}}G_{n}^{\beta }S_{n}^{-1}=V_{n}^{-1}{{C}_{n}}{{\left( {{\left( 1+x \right)}^{n\beta }},x \right)}^{T}}C_{n}^{-1}{{V}_{n}},$$
$$T_{n}^{\beta }={{\tilde{F}}_{n}}{{E}^{n\beta }}\tilde{F}_{n}^{-1}={{\tilde{S}}_{n}}A_{n}^{\beta }\tilde{S}_{n}^{-1}=\tilde{V}_{n}^{-1}{{\tilde{C}}_{n}}\tilde{D}{{\left( {{\left( 1+x \right)}^{n\beta }},x \right)}^{T}}{{\tilde{D}}^{-1}}\tilde{C}_{n}^{-1}{{\tilde{V}}_{n}}.$$
Let $_{\left( \beta  \right)}{{g}_{n}}\left( x \right)$, $_{\left( \beta  \right)}{{h}_{n}}\left( x \right)$ are the numerator polynomials of the matrices 
$$\left( b\left( x{}_{\left( \beta  \right)}{{a}^{\beta }}\left( x \right) \right)\left( 1+x\beta {{\left( \log {}_{\left( \beta  \right)}a\left( x \right) \right)}^{\prime }} \right),x{}_{\left( \beta  \right)}a\left( x \right) \right),$$
$${{\left( b\left( x{}_{\left( \beta  \right)}{{a}^{\beta }}\left( x \right) \right)\left( 1+x\beta {{\left( \log {}_{\left( \beta  \right)}a\left( x \right) \right)}^{\prime }} \right),x{}_{\left( \beta  \right)}a\left( x \right) \right)}_{{{e}^{x}}}},$$
respectively, $_{\left( \beta  \right)}{{\varphi }_{n}}\left( x \right)$ are the numerator polynomials of the matrix ${{\left( 1,x{}_{\left( \beta  \right)}a\left( x \right) \right)}_{{{e}^{x}}}}$. Then
$$G_{n}^{\beta }{{g}_{n}}\left( x \right)={}_{\left( \beta  \right)}{{g}_{n}}\left( x \right),  \quad H_{n}^{\beta }{{h}_{n}}\left( x \right)={}_{\left( \beta  \right)}{{h}_{n}}\left( x \right), \quad T_{n}^{\beta }{{\tilde{\varphi }}_{n}}\left( x \right)={}_{\left( \beta  \right)}{{\tilde{\varphi }}_{n}}\left( x \right).$$
Matrices $G_{n}^{\beta }$, $H_{n}^{\beta }$, $T_{n}^{\beta }$ are characterized by the fact that, in comparison with them, the columns and rows of the matirices $G_{n}^{-\beta }$, $H_{n}^{-\beta }$, $T_{n}^{-\beta }$ are rearranged in the reverse order. Columns of the matrix $G_{n}^{\beta }$ are expressed by the general formula:
$$G_{n}^{\beta }{{x}^{p}}=\sum\limits_{m=0}^{n}{\left( \begin{matrix}
   -n\beta +p  \\
   m  \\
\end{matrix} \right)\left( \begin{matrix}
   n\beta +n-p  \\
   n-m  \\
\end{matrix} \right){{x}^{m}}}.$$

\section{Generalized Euler polynomials}

Let
 $${{u}_{n}}\left( x \right)=\sum\limits_{p=1}^{n}{{{u}_{p}}}{{x}^{p}}, \qquad n>0.$$
Since
$${{a}^{m}}\left( x \right)=\sum\limits_{n=0}^{\infty }{\frac{{{u}_{n}}\left( m \right)}{n!}{{x}^{n}}},  \qquad{{u}_{0}}\left( x \right)=1,$$
then
$$\frac{{{\alpha }_{n}}\left( x \right)}{{{\left( 1-x \right)}^{n+1}}}=\sum\limits_{m=0}^{\infty }{\frac{{{u}_{n}}\left( m \right)}{n!}}{{x}^{m}}=\frac{1}{n!}\sum\limits_{m=0}^{\infty }{{{x}^{m}}}\sum\limits_{p=1}^{n}{{{u}_{p}}}{{m}^{p}}=\frac{1}{n!}\sum\limits_{p=1}^{n}{\sum\limits_{m=0}^{\infty }{{{u}_{p}}{{m}^{p}}{{x}^{m}}}}=$$
$$=\frac{1}{n!}\sum\limits_{p=1}^{n}{\frac{{{u}_{p}}{{A}_{p}}\left( x \right)}{{{\left( 1-x \right)}^{p+1}}}}=\frac{\frac{1}{n!}\sum\limits_{p=1}^{n}{{{u}_{p}}{{\left( 1-x \right)}^{n-p}}{{A}_{p}}\left( x \right)}}{{{\left( 1-x \right)}^{n+1}}}.$$
We introduce the matrices   ${{\tilde{U}}_{n}}$:
$${{\tilde{U}}_{n}}{{x}^{p}}=\frac{1}{n!}{{\left( 1-x \right)}^{n-1-p}}{{\tilde{A}}_{p+1}}\left( x \right),   \qquad p=0, \text{ }1,\text{ } … , \text{ } n-1.$$
For example,
$${{\tilde{U}}_{2}}=\frac{1}{2}\left( \begin{matrix}
   1 & 1  \\
   -1 & 1  \\
\end{matrix} \right), \qquad {{\tilde{U}}_{3}}=\frac{1}{3!}\left( \begin{matrix}
  1 & 1 & 1  \\
   -2 & 0 & 4  \\
  1 & -1 & 1  \\
\end{matrix} \right), \qquad {{\tilde{U}}_{4}}=\frac{1}{4!}\left( \begin{matrix}
   1 & 1 & 1 & 1  \\
   -3 & -1 & 3 & 11  \\
   3 & -1 & -3 & 11  \\
   -1 & 1 & -1 & 1  \\
\end{matrix} \right).$$
Then 
$${{\tilde{U}}_{n}}{{\tilde{u}}_{n}}\left( x \right)={{\tilde{\alpha }}_{n}}\left( x \right).$$
Since 
$$\frac{{{x}^{p+1}}}{{{\left( 1-x \right)}^{n+1}}}=\sum\limits_{m=0}^{\infty }{\left( \begin{matrix}
   m+n-p-1  \\
   n  \\
\end{matrix} \right)}{{x}^{m}}=\sum\limits_{m=0}^{\infty }{\frac{{{\left[ m-p \right]}_{n}}}{n!}}{{x}^{m}}, \qquad 0\le p<n,$$
then
$$\tilde{U}_{n}^{-1}{{x}^{p}}=\frac{1}{x}\prod\limits_{i=0}^{n-1}{\left( x-p+i \right)}={{\left( x-1 \right)}_{p}}{{\left[ x+1 \right]}_{n-p-1}}.$$
For example,
$$\tilde{U}_{2}^{-1}=\left( \begin{matrix}
   1 & -1  \\
   1 & {  }1  \\
\end{matrix} \right),   \qquad\tilde{U}_{3}^{-1}=\left( \begin{matrix}
   2 & -1 & {  }2  \\
   3 & {  }0 & -3  \\
   1 & {  }1 & {  }1  \\
\end{matrix} \right),   \qquad\tilde{U}_{4}^{-1}=\left( \begin{matrix}
   6 & -2 & {  }2 & -6  \\
   11 & -1 & -1 & {  }11  \\
   6 & {  }2 & -2 & -6  \\
   1 & {  }1 & {  }1 & {  }1  \\
\end{matrix} \right).$$
We introduce the matrices ${{J}_{n}}$ corresponding to the operator rearranging the coefficients of the polynomial of degree $n$ in the reverse order. For example,
$${{J}_{3}}=\left( \begin{matrix}
   0 & 0 & 0 & 1  \\
   0 & 0 & 1 & 0  \\
   0 & 1 & 0 & 0  \\
   1 & 0 & 0 & 0  \\
\end{matrix} \right).$$
Denote ${{\tilde{J}}_{n}}={{J}_{n-1}}$.\\
{\bfseries Theorem 1.}
$${{\tilde{U}}_{n}}\left( 1,-x \right)\tilde{U}_{n}^{-1}={{\left( -1 \right)}^{n-1}}{{\tilde{J}}_{n}}.$$
{\bfseries Proof.}
$$\left( 1,-x \right){{\left( x-1 \right)}_{p}}{{\left[ x+1 \right]}_{n-p-1}}={{\left( -x-1 \right)}_{p}}{{\left[ -x+1 \right]}_{n-p-1}}={{\left( -1 \right)}^{n-1}}{{\left( x-1 \right)}_{n-p-1}}{{\left[ x+1 \right]}_{p}},$$
or
$$\left( 1,-x \right)\tilde{U}_{n}^{-1}{{x}^{p}}={{\left( -1 \right)}^{n+1}}\tilde{U}_{n}^{-1}{{x}^{n-p-1}},   \qquad\left( 1,-x \right)\tilde{U}_{n}^{-1}={{\left( -1 \right)}^{n+1}}\tilde{U}_{n}^{-1}{{\tilde{J}}_{n}}.$$
Thus,
$${{\left( -1 \right)}^{n-1}}{{\tilde{J}}_{n}}{{\tilde{\alpha }}_{n}}\left( x \right)={{\tilde{U}}_{n}}{{\tilde{u}}_{n}}\left( -x \right).$$
Denote
$$\left[ n,\searrow  \right]\left( 1,x{{a}^{-1}}\left( x \right) \right)=\frac{\alpha _{n}^{\left( -1 \right)}\left( x \right)}{{{\left( 1-x \right)}^{n+1}}}.$$
Since
$$\left[ n,\to  \right]{{\left( 1,\log {{a}^{-1}}\left( x \right) \right)}_{{{e}^{x}}}}={{u}_{n}}\left( -x \right),$$
then
$$\alpha _{n}^{\left( -1 \right)}\left( x \right)={{\left( -1 \right)}^{n}}x{{J}_{n}}{{\alpha }_{n}}\left( x \right).$$
{\bfseries Theorem 2.}
$${{\alpha }_{n}}\left( 1 \right)={{\left( {{a}_{1}} \right)}^{n}}.$$
{\bfseries Proof.} Denote ${{\tilde{U}}_{n}}{{x}^{p}}={{\tilde{U}}_{p}}\left( x \right)$. Since
$${{a}_{1}}=\left[ x \right]\log a\left( x \right), \quad{{\left( {{a}_{1}} \right)}^{n}}=\left[ {{x}^{n}} \right]{{u}_{n}}\left( x \right);  \quad{{\tilde{U}}_{p}}\left( 1 \right)=0, \quad p<n-1;   \quad{{\tilde{U}}_{n-1}}\left( 1 \right)=1.$$
then
$${{\alpha }_{n}}\left( 1 \right)=\sum\limits_{p=0}^{n-1}{{{u}_{p+1}}{{{\tilde{U}}}_{p}}}\left( 1 \right)={{u}_{n}}={{\left( {{a}_{1}} \right)}^{n}}.$$

The case when ${{a}_{1}}=0$,  the degree of  polynomial ${{u}_{n}}\left( x \right)$ is less than $n$ and matrix ${{\left( 1,\log a\left( x \right) \right)}_{{{e}^{x}}}}$ has no inverse, is possible. This possibility is reflected in the following theorem.\\
{\bfseries Theorem 3.} \emph{
If ${{s}_{n-m}}\left( x \right)$ is a polynomial of degree $n-m-1$, then
$${{\tilde{U}}_{n}}{{s}_{n-m}}\left( x \right)={{\left( 1-x \right)}^{m}}\frac{\left( n-m \right)!}{n!}{{\tilde{U}}_{n-m}}{{s}_{n-m}}\left( x \right).$$
Respectively, if  ${{c}_{n-m}}\left( x \right)$ is a polynomial of degree $<n-m$, then
$$\tilde{U}_{n}^{-1}{{\left( 1-x \right)}^{m}}{{c}_{n-m}}\left( x \right)=\frac{n!}{\left( n-m \right)!}\tilde{U}_{n-m}^{-1}{{c}_{n-m}}\left( x \right).$$}
{\bfseries Proof.} Let ${{I}_{n}}$ is the identity square matrix of order $n+1$, corresponding to the operator annihilating excess columns or rows of matrices. Denote ${{\tilde{I}}_{n}}={{I}_{n-1}}$. Then it is obvious that
$$\left( {{\left( 1-x \right)}^{-m}},x \right){{\tilde{U}}_{n}}{{\tilde{I}}_{n-m}}=\frac{\left( n-m \right)!}{n!}{{\tilde{U}}_{n-m}},$$
or
$${{\tilde{U}}_{n}}{{\tilde{I}}_{n-m}}=\left( {{\left( 1-x \right)}^{m}},x \right)\frac{\left( n-m \right)!}{n!}{{\tilde{U}}_{n-m}}.$$
Respectively,
$$\tilde{U}_{n}^{-1}\left( {{\left( 1-x \right)}^{m}},x \right){{\tilde{I}}_{n-m}}=\frac{n!}{\left( n-m \right)!}\tilde{U}_{n-m}^{-1}.$$

We introduce the matrices   ${{\tilde{V}}_{n}}={{\tilde{J}}_{n}}E{{\tilde{J}}_{n}}=\left( {{\left( 1+x \right)}^{n}},x \right){{P}^{-1}}{{\tilde{I}}_{n}}$. For example,
$${{\tilde{V}}_{4}}=\left( \begin{matrix}
   1 & 0 & 0 & 0  \\
   3 & 1 & 0 & 0  \\
   3 & 2 & 1 & 0  \\
   1 & 1 & 1 & 1  \\
\end{matrix} \right),  \qquad\tilde{V}_{4}^{-1}=\left( \begin{matrix}
   1 & 0 & 0 & 0  \\
   -3 & 1 & 0 & 0  \\
   3 & -2 & 1 & 0  \\
   -1 & 1 & -1 & 1  \\
\end{matrix} \right).$$
Then, as we found in the Introduction,
$$\tilde{V}_{n}^{-1}{{\tilde{v}}_{n}}\left( x \right)={{\tilde{\alpha }}_{n}}\left( x \right),$$
and, hence
   $$\tilde{U}_{n}^{-1}\tilde{V}_{n}^{-1}{{\tilde{v}}_{n}}\left( x \right)={{\tilde{u}}_{n}}\left( x \right).$$
By Theorem 3 we find:
$$\tilde{U}_{n}^{-1}\tilde{V}_{n}^{-1}{{x}^{p}}=\tilde{U}_{n}^{-1}{{\left( 1-x \right)}^{n-p-1}}{{x}^{p}}=\frac{n!}{\left( p+1 \right)!}{{\left( x-1 \right)}_{p}}=$$
$$=\frac{n!}{\left( p+1 \right)!}\sum\limits_{m=1}^{p+1}{s\left( p+1,\text{ }m \right){{x}^{m-1}}},$$
where $s\left( p+1,\text{ }m \right)$ are the Stirling numbers of  the first kind. Hence
$${{\tilde{V}}_{n}}{{\tilde{U}}_{n}}{{x}^{p}}=\frac{1}{n!}\sum\limits_{m=1}^{p+1}{m!S\left( p+1,\text{ }m \right)}\text{ }{{x}^{m-1}},$$
where $S\left( p+1,\text{ }m \right)$ are the Stirling numbers of the second kind. For example,
$$\tilde{U}_{4}^{-1}\tilde{V}_{4}^{-1}=4!\left( \begin{matrix}
   1 & -1 & {  }2 & -6  \\
   0 & {  }1 & -3 & { }11  \\
   0 & {  }0 & {  }1 & -6  \\
   0 & {  }0 & {  }0 & {  }1  \\
\end{matrix} \right)\left( \begin{matrix}
   1 & 0 & 0 & 0  \\
   0 & \frac{1}{2} & 0 & 0  \\
   0 & 0 & \frac{1}{3!} & 0  \\
   0 & 0 & 0 & \frac{1}{4!}  \\
\end{matrix} \right),$$
$${{\tilde{V}}_{4}}{{\tilde{U}}_{4}}=\frac{1}{4!}\left( \begin{matrix}
   1 & 0 & 0 & 0  \\
   0 & 2 & 0 & 0  \\
   0 & 0 & 3! & 0  \\
   0 & 0 & 0 & 4!  \\
\end{matrix} \right)\left( \begin{matrix}
   1 & 1 & 1 & 1  \\
   0 & 1 & 3 & 7  \\
   0 & 0 & 1 & 6  \\
   0 & 0 & 0 & 1  \\
\end{matrix} \right).$$

Generalized binomial series,
$$_{\left( \beta  \right)}{{a}^{\varphi }}\left( x \right)=\sum\limits_{n=0}^{\infty }{\frac{\varphi }{\varphi +n\beta }}\left( \begin{matrix}
   \varphi +n\beta   \\
   n  \\
\end{matrix} \right){{x}^{n}};$$
$$_{\left( 0 \right)}a\left( x \right)=1+x,  \qquad_{\left( 1 \right)}a\left( x \right)={{\left( 1-x \right)}^{-1}},    \qquad_{\left( 2 \right)}a\left( x \right)=\frac{1-{{\left( 1-4x \right)}^{{1}/{2}\;}}}{2x},$$
$$_{\left( -1 \right)}a\left( x \right)=\frac{1+{{\left( 1+4x \right)}^{{1}/{2}\;}}}{2},    \qquad_{\left( {1}/{2}\; \right)}a\left( x \right)={{\left( \frac{x}{2}+{{\left( 1+\frac{{{x}^{2}}}{4} \right)}^{{1}/{2}\;}} \right)}^{2}},$$
takes important place in our studies. Generalization that underlies it can be extended to each formal power series $a\left( x \right)$, ${{a}_{0}}=1$. Each such series  is associated by means of the Lagrange transform
$${{a}^{\varphi }}\left( x \right)=\sum\limits_{n=0}^{\infty }{\frac{{{x}^{n}}}{{{a}^{\beta n}}\left( x \right)}\left[ {{x}^{n}} \right]}\left( 1-x\beta {{\left( \log a\left( x \right) \right)}^{\prime }} \right){{a}^{\varphi +\beta n}}\left( x \right)$$
with the set of series $_{\left( \beta  \right)}a\left( x \right)$, $_{\left( 0 \right)}a\left( x \right)=a\left( x \right)$,   such that
$${}_{\left( \beta  \right)}a\left( x{{a}^{-\beta }}\left( x \right) \right)=a\left( x \right),   \qquad a\left( x{}_{\left( \beta  \right)}{{a}^{\beta }}\left( x \right) \right)={}_{\left( \beta  \right)}a\left( x \right),$$
$$\left[ {{x}^{n}} \right]{}_{\left( \beta  \right)}{{a}^{\varphi }}\left( x \right)=\left[ {{x}^{n}} \right]\left( 1-x\beta \frac{{a}'\left( x \right)}{a\left( x \right)} \right){{a}^{\varphi +\beta n}}\left( x \right)=\frac{\varphi }{\varphi +\beta n}\left[ {{x}^{n}} \right]{{a}^{\varphi +\beta n}}\left( x \right),$$
$$\left[ {{x}^{n}} \right]\left( 1+x\beta \frac{_{\left( \beta  \right)}{a}'\left( x \right)}{_{\left( \beta  \right)}a\left( x \right)} \right){}_{\left( \beta  \right)}{{a}^{\varphi }}\left( x \right)=\frac{\varphi +\beta n}{\varphi }\left[ {{x}^{n}} \right]{}_{\left( \beta  \right)}{{a}^{\varphi }}\left( x \right)=\left[ {{x}^{n}} \right]a_{{}}^{\varphi +\beta n}\left( x \right),$$
$${{\left( 1,x{}_{\left( \beta  \right)}{{a}^{\varphi }}\left( x \right) \right)}^{-1}}=\left( 1,{{x}_{\left( \beta -\varphi  \right)}}{{a}^{-\varphi }}\left( x \right) \right),$$
$${{\left( 1+x\varphi {{\left( \log {}_{\left( \beta  \right)}a\left( x \right) \right)}^{\prime }},x{}_{\left( \beta  \right)}{{a}^{\varphi }}\left( x \right) \right)}^{-1}}=\left( 1-x\varphi {{\left( \log {}_{\left( \beta -\varphi  \right)}a\left( x \right) \right)}^{\prime }},x{}_{\left( \beta -\varphi  \right)}{{a}^{-\varphi }}\left( x \right) \right).$$

Denote
$$\left[ n,\to  \right]{{\left( 1,\log {}_{\left( \beta  \right)}a\left( x \right) \right)}_{{{e}^{x}}}}={}_{\left( \beta  \right)}{{u}_{n}}\left( x \right).$$
Then 
$$_{\left( \beta  \right)}{{a}^{\varphi }}\left( x \right)=\sum\limits_{n=0}^{\infty }{\frac{\varphi }{\varphi +n\beta }}\frac{{{u}_{n}}\left( \varphi +n\beta  \right)}{n!}{{x}^{n}},$$
$$_{\left( \beta  \right)}{{u}_{n}}\left( x \right)=x{{\left( x+n\beta  \right)}^{-1}}{{u}_{n}}\left( x+n\beta  \right).$$
Let 
$${{\left( 1,\log a\left( x \right) \right)}^{-1}}=\left( 1,q\left( x \right) \right).$$
Then
$${{\left( 1,\log {}_{\left( \beta  \right)}a\left( x \right) \right)}^{-1}}=\left( 1,q\left( x \right){{e}^{-\beta x}} \right).$$
Denote
$${{\left( 1,q\left( x \right){{e}^{-\beta x}} \right)}_{{{e}^{x}}}}{{x}^{n}}={}_{\left( \beta  \right)}{{q}_{n}}\left( x \right).$$
Then
$$_{\left( \beta  \right)}{{q}_{n}}\left( x \right)={{\left( 1+n\beta x \right)}^{-1}}{{q}_{n}}\left( \frac{x}{1+n\beta x} \right),$$
$$\sum\limits_{n=0}^{\infty }{_{\left( \beta  \right)}{{u}_{n}}}\left( \varphi  \right){}_{\left( \beta  \right)}{{q}_{n}}\left( x \right)={{\left( 1-\varphi x \right)}^{-1}}.$$

Series $_{\left( \beta  \right)}a\left( x \right)$ for integer $\beta $, denoted by ${{S}_{\beta }}\left( x \right)$, were introduced in [15]. In [16] these series, called generalized Lagrange series, are considered in connection with the Riordan arrays. Properties of these series intersect with the properties of Sheffer sequences, therefore the identities associated with them  can be found in the umbral calculus.

Denote
$$\frac{1}{x}{}_{\left( \beta  \right)}{{u}_{n}}\left( x \right)={}_{\left( \beta  \right)}{{\tilde{u}}_{n}}\left( x \right),
\quad\left[ n,\searrow  \right]\left( 1,x{}_{\left( \beta  \right)}a\left( x \right) \right)=\frac{_{\left( \beta  \right)}{{\alpha }_{n}}\left( x \right)}{{{\left( 1-x \right)}^{n+1}}},  \quad\frac{1}{x}{}_{\left( \beta  \right)}{{\alpha }_{n}}\left( x \right)={}_{\left( \beta  \right)}{{\tilde{\alpha }}_{n}}\left( x \right).$$
Then
$${{E}^{n\beta }}{{\tilde{u}}_{n}}\left( x \right)={{\tilde{u}}_{n}}\left( x+n\beta  \right)={}_{\left( \beta  \right)}{{\tilde{u}}_{n}}\left( x \right),
\qquad{{U}_{n}}{{E}^{n\beta }}U_{n}^{-1}{{\tilde{\alpha }}_{n}}\left( x \right)={}_{\left( \beta  \right)}{{\tilde{\alpha }}_{n}}\left( x \right).$$

Denote
$${{\tilde{U}}_{n}}{{E}^{n\beta }}\tilde{U}_{n}^{-1}=A_{n}^{\beta }.$$
Since
$$\left( 1,-x \right){{E}^{n\beta }}\left( 1,-x \right)={{E}^{-n\beta }},   \qquad{{\tilde{U}}_{n}}\left( 1,-x \right)\tilde{U}_{n}^{-1}={{\left( -1 \right)}^{n-1}}{{\tilde{J}}_{n}},$$
then
$${{\tilde{J}}_{n}}A_{n}^{\beta }{{\tilde{J}}_{n}}=A_{n}^{-\beta }.$$
For example,
$${{A}_{2}}=\left( \begin{matrix}
  2 & 1  \\
   -1 & 0  \\
\end{matrix} \right),   \quad{{A}_{3}}=\left( \begin{matrix}
 5 & {5}/{2}\; & 1  \\
   -6 & -2 & 0  \\
  2 & {1}/{2}\; & 0  \\
\end{matrix} \right),    
\quad{{A}_{4}}=\left( \begin{matrix}
14 & 7 & 3 & 1  \\
   -28 & -{35}/{3}\; & -{10}/{3}\; & 0  \\
  20 & {22}/{3}\; & {5}/{3}\; & 0  \\
   -5 & -{5}/{3}\; & -{1}/{3}\; & 0  \\
\end{matrix} \right);$$
$$A_{2}^{-1}=\left( \begin{matrix}
   0 & -1  \\
   1 & 2  \\
\end{matrix} \right),   \quad A_{3}^{-1}=\left( \begin{matrix}
   0 & {1}/{2}\; & 2  \\
   0 & -2 & -6  \\
   1 & {5}/{2}\; & 5  \\
\end{matrix} \right),   \quad A_{4}^{-1}=\left( \begin{matrix}
   0 & -{1}/{3}\; & -{5}/{3}\; & -5  \\
   0 & {5}/{3}\; & {22}/{3}\; & 20  \\
   0 & -{10}/{3}\; & -{35}/{3}\; & -28  \\
   1 & 3 & 7 & 14  \\
\end{matrix} \right).$$
Since  $\left[ x \right]{}_{\left( \beta  \right)}a\left( x \right)=\left[ x \right]a\left( x \right)={{a}_{1}}$, then ${}_{\left( \beta  \right)}{{\alpha }_{n}}\left( 1 \right)={{\alpha }_{n}}\left( 1 \right)$ and sum of the elements of each column of the matrix $A_{n}^{\beta }$ is $1$. From Theorem 3 it follows that
$$\left( {{\left( 1-x \right)}^{-m}},x \right)A_{n}^{\beta }\left( {{\left( 1-x \right)}^{m}},x \right){{\tilde{I}}_{n-m}}=A_{n-m}^{\frac{n\beta }{n-m}}.$$

We introduce the diagonal matrix $\tilde{D}$, $\tilde{D}{{x}^{n}}=\left( n+1 \right){{x}^{n}}$:
$$\tilde{D}=\left( \begin{matrix}
   1 & 0 & 0 & \cdots   \\
   0 & 2 & 0 & \cdots   \\
   0 & 0 & 3 & \cdots   \\
   \vdots  & \vdots  & \vdots  & \ddots   \\
\end{matrix} \right),   \qquad{{\tilde{D}}^{-1}}=\left( \begin{matrix}
   1 & 0 & 0 & \cdots   \\
   0 & \frac{1}{2} & 0 & \cdots   \\
   0 & 0 & \frac{1}{3} & \cdots   \\
   \vdots  & \vdots  & \vdots  & \ddots   \\
\end{matrix} \right).$$
{\bfseries Theorem 4.}
$$A_{n}^{\beta }=\tilde{V}_{n}^{-1}\tilde{D}{{\left( {{\left( 1+x \right)}^{n\beta }},x \right)}^{T}}{{\tilde{D}}^{-1}}{{\tilde{V}}_{n}}.$$
{\bfseries Proof.} Columns of the matrices ${{\tilde{V}}_{n}}{{\tilde{U}}_{n}}$, $\tilde{U}_{n}^{-1}\tilde{V}_{n}^{-1}$ are connected a certain way with the rows of the matrices ${{\left( {{e}^{x}},{{e}^{x}}-1 \right)}_{{{e}^{x}}}}$, ${{\left( {{\left( 1+x \right)}^{-1}},\log \left( 1+x \right) \right)}_{{{e}^{x}}}}$:
$$n!\left| {{e}^{x}} \right|{{\tilde{D}}^{-1}}{{\tilde{V}}_{n}}{{\tilde{U}}_{n}}{{x}^{p}}=\left[ p,\to  \right]{{\left( {{e}^{x}},{{e}^{x}}-1 \right)}_{{{e}^{x}}}},$$
$$\left( {1}/{n!}\; \right)\tilde{U}_{n}^{-1}\tilde{V}_{n}^{-1}{{\left| {{e}^{x}} \right|}^{-1}}\tilde{D}{{x}^{p}}=\left[ p,\to  \right]{{\left( {{\left( 1+x \right)}^{-1}},\log \left( 1+x \right) \right)}_{{{e}^{x}}}}.$$
Since
$${{\left( {{\left( 1+x \right)}^{-1}},\log \left( 1+x \right) \right)}_{{{e}^{x}}}}{{\left( {{e}^{n\beta }},x \right)}_{{{e}^{x}}}}{{\left( {{e}^{x}},{{e}^{x}}-1 \right)}_{{{e}^{x}}}}={{\left( {{\left( 1+x \right)}^{n\beta }},x \right)}_{{{e}^{x}}}},$$
then
$${{\tilde{V}}_{n}}{{\tilde{U}}_{n}}{{E}^{n\beta }}\tilde{U}_{n}^{-1}\tilde{V}_{n}^{-1}=\tilde{D}{{\left( {{\left( 1+x \right)}^{n\beta }},x \right)}^{T}}{{\tilde{D}}^{-1}}{{\tilde{I}}_{n}}.$$

 Thus,
$${{A}_{2}}=\left( \begin{matrix}
   1 & 0  \\
   -1 & 1  \\
\end{matrix} \right)\left( \begin{matrix}
   1 & 0  \\
   0 & 2  \\
\end{matrix} \right)\left( \begin{matrix}
   1 & 2  \\
   0 & 1  \\
\end{matrix} \right)\left( \begin{matrix}
   1 & 0  \\
   0 & \frac{1}{2}  \\
\end{matrix} \right)\left( \begin{matrix}
   1 & 0  \\
   1 & 1  \\
\end{matrix} \right),$$
$${{A}_{3}}=\left( \begin{matrix}
   1 & 0 & 0  \\
   -2 & 1 & 0  \\
   1 & -1 & 1  \\
\end{matrix} \right)\left( \begin{matrix}
   1 & 0 & 0  \\
   0 & 2 & 0  \\
   0 & 0 & 3  \\
\end{matrix} \right)\left( \begin{matrix}
   1 & 3 & 3  \\
   0 & 1 & 3  \\
   0 & 0 & 1  \\
\end{matrix} \right)\left( \begin{matrix}
   1 & 0 & 0  \\
   0 & \frac{1}{2} & 0  \\
   0 & 0 & \frac{1}{3}  \\
\end{matrix} \right)\left( \begin{matrix}
   1 & 0 & 0  \\
   2 & 1 & 0  \\
   1 & 1 & 1  \\
\end{matrix} \right),$$
$${{A}_{4}}=\left( \begin{matrix}
   1 & 0 & 0 & 0  \\
   -3 & 1 & 0 & 0  \\
   3 & -2 & 1 & 0  \\
   -1 & 1 & -1 & 1  \\
\end{matrix} \right)\left( \begin{matrix}
   1 & 0 & 0 & 0  \\
   0 & 2 & 0 & 0  \\
   0 & 0 & 3 & 0  \\
   0 & 0 & 0 & 4  \\
\end{matrix} \right)\left( \begin{matrix}
   1 & 4 & 6 & 4  \\
   0 & 1 & 4 & 6  \\
   0 & 0 & 1 & 4  \\
   0 & 0 & 0 & 1  \\
\end{matrix} \right)\left( \begin{matrix}
   1 & 0 & 0 & 0  \\
   0 & \frac{1}{2} & 0 & 0  \\
   0 & 0 & \frac{1}{3} & 0  \\
   0 & 0 & 0 & \frac{1}{4}  \\
\end{matrix} \right)\left( \begin{matrix}
   1 & 0 & 0 & 0  \\
   3 & 1 & 0 & 0  \\
   3 & 2 & 1 & 0  \\
   1 & 1 & 1 & 1  \\
\end{matrix} \right).$$

Denote
$$\left[ n,\to  \right]\left( 1,{}_{\left( \beta  \right)}a\left( x \right)-1 \right)={}_{\left( \beta  \right)}{{v}_{n}}\left( x \right), \qquad\frac{1}{x}{}_{\left( \beta  \right)}{{v}_{n}}\left( x \right)={}_{\left( \beta  \right)}{{\tilde{v}}_{n}}\left( x \right).$$
Then
$$\tilde{D}{{\left( {{\left( 1+x \right)}^{n\beta }},x \right)}^{T}}{{\tilde{D}}^{-1}}{{\tilde{v}}_{n}}\left( x \right)={}_{\left( \beta  \right)}{{\tilde{v}}_{n}}\left( x \right).$$

Let $_{\left( \beta  \right)}a\left( x \right)$ is the generalized binomial series. Then $a\left( x \right)=1+x$,
$$\frac{_{\left( \beta  \right)}{{\alpha }_{n}}\left( x \right)}{{{\left( 1-x \right)}^{n+1}}}=\sum\limits_{m=0}^{\infty }{\frac{m}{m+n\beta }}\left( \begin{matrix}
   m+n\beta   \\
   n  \\
\end{matrix} \right){{x}^{m}}.$$
{\bfseries Theorem 5.}
$$_{\left( \beta  \right)}{{\alpha }_{n}}\left( x \right)=\frac{1}{n}\sum\limits_{m=1}^{n}{\left( \begin{matrix}
   n\left( 1-\beta  \right)  \\
   m-1  \\
\end{matrix} \right)\left( \begin{matrix}
   n\beta   \\
   n-m  \\
\end{matrix} \right){{x}^{m}}}.$$
{\bfseries Proof.} We use the factorial representation of binomial coefficients, i.e. we prove the theorem for positive integers $\beta $. By polynomial argument (binomial coefficients under consideration are polynomials in $\beta $) this is equivalent to the general proof . Since ${{\tilde{v}}_{n}}\left( x \right)={{x}^{n-1}}$, then
$$_{\left( \beta  \right)}{{\tilde{v}}_{n}}\left( x \right)=\tilde{D}{{\left( {{\left( 1+x \right)}^{n\beta }},x \right)}^{T}}{{\tilde{D}}^{-1}}{{x}^{n-1}}=\sum\limits_{m=0}^{n-1}{\frac{m+1}{n}\left( \begin{matrix}
   n\beta   \\
   n-m-1  \\
\end{matrix} \right)}{{x}^{m}}.$$
Since
$$\left[ m,\to  \right]\tilde{V}_{n}^{-1}=\sum\limits_{i=0}^{m}{\left( \begin{matrix}
   m-n  \\
   m-i  \\
\end{matrix} \right){{x}^{i}}}=\sum\limits_{i=0}^{m}{{{\left( -1 \right)}^{m-i}}}\left( \begin{matrix}
   n-i-1  \\
   m-i  \\
\end{matrix} \right){{x}^{i}},$$
then
$$\left[ {{x}^{m}} \right]{}_{\left( \beta  \right)}{{\tilde{\alpha }}_{n}}\left( x \right)=\left[ {{x}^{m}} \right]\tilde{V}{{_{n}^{-1}}_{\left( \beta  \right)}}{{\tilde{v}}_{n}}\left( x \right)=$$
$$=\sum\limits_{i=0}^{m}{{{\left( -1 \right)}^{m-i}}}\left( \begin{matrix}
   n-i-1  \\
   m-i  \\
\end{matrix} \right)\frac{\left( i+1 \right)}{n}\left( \begin{matrix}
   n\beta   \\
   n-i-1  \\
\end{matrix} \right)\frac{\left( n\beta -n+m+1 \right)!}{\left( n\beta -n+m+1 \right)!}=$$
$$=\frac{1}{n}\left( \begin{matrix}
   n\beta   \\
   n-m-1  \\
\end{matrix} \right)\sum\limits_{i=0}^{m}{{{\left( -1 \right)}^{m-i}}\left( i+1 \right)}\left( \begin{matrix}
   n\beta -n+m+1  \\
   m-i  \\
\end{matrix} \right)=$$ 
$$=\frac{1}{n}\left( \begin{matrix}
   n\beta   \\
   n-m-1  \\
\end{matrix} \right){{\left( -1 \right)}^{m}}\left( \begin{matrix}
   n\beta -n+m-1  \\
   m  \\
\end{matrix} \right)=\frac{1}{n}\left( \begin{matrix}
   n\beta   \\
   n-m-1  \\
\end{matrix} \right)\left( \begin{matrix}
   n\left( 1-\beta  \right)  \\
   m  \\
\end{matrix} \right).$$

Note that
$$_{\left( 0 \right)}{{\alpha }_{n}}\left( x \right)={{x}^{n}}, \qquad_{\left( 1 \right)}{{\alpha }_{n}}\left( x \right)=x,  \qquad_{\left( {1}/{2}\; \right)}{{\alpha }_{2n}}\left( x \right)=\frac{1}{2}\left( 1+x \right){{x}^{n}};$$
since
$$_{\left( 1-\beta  \right)}a\left( x \right)={}_{\left( \beta  \right)}{{a}^{-1}}\left( -x \right),$$  
then
$$_{\left( 1-\beta  \right)}{{\alpha }_{n}}\left( x \right)=x{{J}_{n}}{}_{\left( \beta  \right)}{{\alpha }_{n}}\left( x \right).$$
\section{Generalized Narayana polynomials}

Constructive relationships between the ordinary and the exponential Riordan arrays exist. Particular manifestations of these relationships resemble the details of construction, the general plan of which is a secret for us. Following [17] – [19], we will consider some of such manifestations associated with the numerator polynomials. Since
$$[n,\searrow ]\left( 1,xa\left( x \right) \right)=\sum\limits_{m=0}^{\infty }{\frac{{{u}_{n}}\left( m \right)}{n!}}{{x}^{m}},$$ 
then
$$[n,\searrow ]{{\left( 1,xa\left( x \right) \right)}_{{{e}^{x}}}}=\sum\limits_{m=0}^{\infty }{\frac{\left( m+n \right)!}{m!}\frac{{{u}_{n}}\left( m \right)}{n!}}{{x}^{m}}=\sum\limits_{m=0}^{\infty }{\frac{{{\left[ m+1 \right]}_{n}}{{u}_{n}}\left( m \right)}{n!}{{x}^{m}}}.$$  
 If ${{a}_{1}}\ne 0$, then ${{\left[ x+1 \right]}_{n}}{{u}_{n}}\left( x \right)$ is the polynomial of degree $2n$, so that
$$\sum\limits_{m=0}^{\infty }{\frac{{{\left[ m+1 \right]}_{n}}{{u}_{n}}\left( m \right)}{n!}{{x}^{m}}}=\frac{{{\varphi }_{n}}\left( x \right)}{{{\left( 1-x \right)}^{2n+1}}},$$
where 
$${{\varphi }_{n}}\left( x \right)=x\frac{\left( 2n \right)!}{n!}{{\tilde{U}}_{2n}}{{\left[ x+1 \right]}_{n}}{{\tilde{u}}_{n}}\left( x \right).$$
Since ${{\left[ x+1 \right]}_{n}}{{u}_{n}}\left( x \right)=0$ when $x=0$, $-1$, … , $-n$, then, in accordance with the Theorem1,
$$\sum\limits_{m=0}^{\infty }{\frac{{{\left[ -m+1 \right]}_{n}}{{u}_{n}}\left( -m \right)}{n!}{{x}^{m}}}=\frac{{{\left( -1 \right)}^{2n}}x{{J}_{2n}}{{\varphi }_{n}}\left( x \right)}{{{\left( 1-x \right)}^{2n+1}}}=\frac{{{\left( -1 \right)}^{2n}}{{x}^{n+1}}{{J}_{n}}{{\varphi }_{n}}\left( x \right)}{{{\left( 1-x \right)}^{2n+1}}},$$
i.e. ${{\varphi }_{n}}\left( x \right)$  is the polynomial of degree $\le n$. Since
$$\left[ {{x}^{2n}} \right]{{\left[ x+1 \right]}_{n}}{{u}_{n}}\left( x \right)={{\left( {{a}_{1}} \right)}^{n}},$$
 then, in accordance with the Theorem 2,
$${{\varphi }_{n}}\left( 1 \right)={{\left( {{a}_{1}} \right)}^{n}}\frac{\left( 2n \right)!}{n!}.$$
If $a\left( x \right)={{\left( 1-x \right)}^{-1}}$, then 
$${{\varphi }_{n}}\left( x \right)=\left( n+1 \right)!{{N}_{n}}\left( x \right)={{\left( 1-x \right)}^{2n+1}}\sum\limits_{m=0}^{\infty }{{{\left[ m+1 \right]}_{n}}\left( \begin{matrix}
   m+n-1  \\
   n  \\
\end{matrix} \right)}{{x}^{m}},$$
where 
$${{N}_{n}}\left( x \right)=\frac{1}{n}\sum\limits_{m=1}^{n}{\left( \begin{matrix}
   n  \\
   m-1  \\
\end{matrix} \right)\left( \begin{matrix}
   n  \\
   n-m  \\
\end{matrix} \right){{x}^{m}}}$$ 
is the Narayana polynomials. In this connection we will called polynomials ${{\varphi }_{n}}\left( x \right)$ the generalized Narayana polynomials (GNP).

Since
$$\frac{{{\varphi }_{n}}\left( t \right)}{\left( n+1 \right)!{{\left( 1-t \right)}^{2n+1}}}=\sum\limits_{m=0}^{\infty }{\frac{1}{n+1}}\left( \begin{matrix}
   n+m  \\
   m  \\
\end{matrix} \right)\left[ {{x}^{n}} \right]{{a}^{m}}\left( x \right){{t}^{m}}=$$
$$=\frac{1}{n+1}\left[ {{x}^{n}} \right]{{\left( 1-ta\left( x \right) \right)}^{-n-1}}=\left[ {{x}^{n}} \right]b\left( x \right),$$
where
$$b\left( x \right)=\frac{1}{1-ta\left( xb\left( x \right) \right)} , \qquad\left( 1,xb\left( x \right) \right)={{\left( 1,x\left( 1-ta\left( x \right) \right) \right)}^{-1}},$$
then
$$\sum\limits_{n=0}^{\infty }{{{\varphi }_{n}}}\left( t \right)\frac{{{x}^{n}}}{\left( n+1 \right)!}=\left( 1-t \right)b\left( x{{\left( 1-t \right)}^{2}} \right).$$
For example,
$$a\left( x \right)=\frac{1}{1-x},   \qquad b\left( x \right)=\frac{1+x-t-\sqrt{1-2x-2t-2xt+{{x}^{2}}+{{t}^{2}}}}{2x},$$
$$\sum\limits_{n=0}^{\infty }{{{\varphi }_{n}}\left( t \right)}\frac{{{x}^{n}}}{\left( n+1 \right)!}=\sum\limits_{n=0}^{\infty }{{{N}_{n}}}\left( t \right){{x}^{n}}=\frac{1+x\left( 1-t \right)-\sqrt{1-2x\left( 1+t \right)+{{x}^{2}}{{\left( 1-t \right)}^{2}}}}{2x}.$$

Generating functions of the numerator polynomials of the matrices ${{\left( 1,{{e}^{x}}-1 \right)}_{{{e}^{x}}}}$, ${{\left( {{e}^{x}},{{e}^{x}}-1 \right)}_{{{e}^{x}}}}$  are considered in [20].

We introduce the matrices ${{\tilde{F}}_{n}}$:
$${{\tilde{F}}_{n}}=\frac{\left( 2n \right)!}{n!}{{\tilde{U}}_{2n}}\left( {{\left[ x+1 \right]}_{n}},x \right){{\tilde{I}}_{n}},  \qquad\tilde{F}_{n}^{-1}=\frac{n!}{\left( 2n \right)!}{{\left( {{\left[ x+1 \right]}_{n}},x \right)}^{-1}}\tilde{U}_{2n}^{-1}{{\tilde{I}}_{n}};$$
$${{\tilde{F}}_{n}}{{x}^{p}}=\frac{\left( 2n \right)!}{n!}{{\tilde{U}}_{2n}}{{x}^{p}}{{\left[ x+1 \right]}_{n}}={{\left( 1-x \right)}^{2n+1}}\sum\limits_{m=1}^{\infty }{{{m}^{p+1}}}\left( \begin{matrix}
   m+n  \\
   n  \\
\end{matrix} \right){{x}^{m-1}},$$
$$\tilde{F}_{n}^{-1}{{x}^{p}}=\frac{n!}{\left( 2n \right)!}{{\left( x-1 \right)}_{p}}{{\left[ x+n+1 \right]}_{n-p-1}},  \qquad p=0, \text{ }1,\text{ } … , \text{ } n-1.$$
For example,
$${{\tilde{F}}_{2}}=3\left( \begin{matrix}
   1 & 1  \\
   -1 & 3  \\
\end{matrix} \right), \qquad{{\tilde{F}}_{3}}=4\left( \begin{matrix}
   1 & 1 & 1  \\
   -2 & 3 & 13  \\
   1 & -4 & 16  \\
\end{matrix} \right), \qquad{{\tilde{F}}_{4}}=5\left( \begin{matrix}
   1 & 1 & 1 & 1  \\
   -3 & 3 & 15 & 39  \\
   3 & -9 & 9 & 171  \\
   -1 & 5 & -25 & 125  \\
\end{matrix} \right);$$
$$\tilde{F}_{2}^{-1}=\frac{2!}{4!}\left( \begin{matrix}
   3 & -1  \\
   1 & 1  \\
\end{matrix} \right),  \quad\tilde{F}_{3}^{-1}=\frac{3!}{6!}\left( \begin{matrix}
   20 & -4 & 2  \\
   9 & 3 & -3  \\
   1 & 1 & 1  \\
\end{matrix} \right),
\quad\tilde{F}_{4}^{-1}=\frac{4!}{8!}\left( \begin{matrix}
   210 & -30 & 10 & -6  \\
   107 & 19 & -13 & 11  \\
   18 & 10 & 2 & -6  \\
   1 & 1 & 1 & 1  \\
\end{matrix} \right).$$
Denote $\left( {1}/{x}\; \right){{\varphi }_{n}}\left( x \right)={{\tilde{\varphi }}_{n}}\left( x \right)$. Then
$${{\tilde{F}}_{n}}{{\tilde{u}}_{n}}\left( x \right)={{\tilde{\varphi }}_{n}}\left( x \right).$$
Let ${{\alpha }_{n}}\left( x \right)|a\left( x \right)$, ${{\varphi }_{n}}\left( x \right)|a\left( x \right)$ denotes respectively GEP, GNP associated with the matrices $\left( 1,xa\left( x \right) \right)$, ${{\left( 1,xa\left( x \right) \right)}_{{{e}^{x}}}}$. Then
$$x{{\tilde{F}}_{n}}{{x}^{n-1}}={{\varphi }_{n}}\left( x \right)|{{e}^{x}}.$$
{\bfseries Theorem 6.}
$${{\tilde{F}}_{n}}{{E}^{n}}\left( 1,-x \right)\tilde{F}_{n}^{-1}={{\left( -1 \right)}^{n-1}}{{\tilde{J}}_{n}}.$$
{\bfseries Proof.} 
$${{E}^{n}}\left( 1,-x \right){{\left( x-1 \right)}_{p}}{{\left[ x+n+1 \right]}_{n-p-1}}={{\left( -x-n-1 \right)}_{p}}{{\left[ -x+1 \right]}_{n-p-1}}=$$
$$={{\left( -1 \right)}^{n-1}}{{\left( x-1 \right)}_{n-p-1}}{{\left[ x+n+1 \right]}_{p}},$$
or
$${{E}^{n}}\left( 1,-x \right)\tilde{F}_{n}^{-1}{{x}^{p}}={{\left( -1 \right)}^{n-1}}\tilde{F}_{n}^{-1}{{x}^{n-1-p}},   \qquad{{E}^{n}}\left( 1,-x \right)\tilde{F}_{n}^{-1}={{\left( -1 \right)}^{n-1}}\tilde{F}_{n}^{-1}{{\tilde{J}}_{n}}.$$

Thus,
$${{\left( -1 \right)}^{n-1}}{{\tilde{J}}_{n}}{{\tilde{\varphi }}_{n}}\left( x \right)={{\tilde{F}}_{n}}{{\tilde{u}}_{n}}\left( -x-n \right).$$
We denote by using the notation for the series $_{\left( \beta  \right)}a\left( x \right)$:
$${{\left( 1,xa\left( x \right) \right)}^{-1}}=\left( 1,x{}_{\left( -1 \right)}{{a}^{-1}}\left( x \right) \right).$$
Then
$$\left[ n,\to  \right]{{\left( 1,\log {}_{\left( -1 \right)}{{a}^{-1}}\left( x \right) \right)}_{{{e}^{x}}}}=-x{{\tilde{u}}_{n}}\left( -x-n \right).$$
Denote
$$[n,\searrow ]\left( 1,xa\left( x \right) \right)_{{{e}^{x}}}^{-1}=\frac{\varphi _{n}^{\left[ -1 \right]}\left( x \right)}{{{\left( 1-x \right)}^{2n+1}}}.$$
Then
$$\varphi _{n}^{\left[ -1 \right]}\left( x \right)={{\left( -1 \right)}^{n}}x{{J}_{n}}{{\varphi }_{n}}\left( x \right).$$

Thus, if the matrix $\left( 1,xa\left( x \right) \right)$ is a pseudo-involution, i.e. ${{\left( 1,xa\left( x \right) \right)}^{-1}}=\left( 1,xa\left( -x \right) \right)$, then  ${{\varphi }_{n}}\left( x \right)=x{{J}_{n}}{{\varphi }_{n}}\left( x \right)$.

We introduce the matrices ${{\tilde{S}}_{n}}$:
$${{\tilde{S}}_{n}}={{\tilde{F}}_{n}}\tilde{U}_{n}^{-1},  \qquad\tilde{S}_{n}^{-1}={{\tilde{U}}_{n}}\tilde{F}_{n}^{-1}.$$
For example,
$${{\tilde{S}}_{2}}=3!\left( \begin{matrix}
   1 & 0  \\
   1 & 2  \\
\end{matrix} \right),  \qquad{{\tilde{S}}_{3}}=4!\left( \begin{matrix}
   1 & 0 & 0  \\
   3 & {5}/{2}\; & 0  \\
   1 & {5}/{2}\; & 5  \\
\end{matrix} \right),  \qquad{{\tilde{S}}_{4}}=5!\left( \begin{matrix}
   1 & 0 & 0 & 0  \\
   6 & 3 & 0 & 0  \\
   6 & 8 & 7 & 0  \\
   1 & 3 & 7 & 14  \\
\end{matrix} \right);$$
$$\tilde{S}_{2}^{-1}=\frac{2!}{4!}\left( \begin{matrix}
   2 & 0  \\
   -1 & 1  \\
\end{matrix} \right), \quad\tilde{S}_{3}^{-1}=\frac{3!}{6!}\left( \begin{matrix}
   5 & 0 & 0  \\
   -6 & 2 & 0  \\
   2 & -1 & 1  \\
\end{matrix} \right), \quad\tilde{S}_{4}^{-1}=\frac{4!}{8!}\left( \begin{matrix}
   14 & 0 & 0 & 0  \\
   -28 & {14}/{3}\; & 0 & 0  \\
   20 & {-16}/{3}\; & 2 & 0  \\
   -5 & {5}/{3}\; & -1 & 1  \\
\end{matrix} \right).$$
Then
$${{\tilde{S}}_{n}}{{\tilde{\alpha }}_{n}}\left( x \right)={{\tilde{\varphi }}_{n}}\left( x \right).$$
{\bfseries Theorem 7.}
$${{\tilde{S}}_{n}}=\tilde{V}_{n}^{-1}{{\tilde{C}}_{n}}{{\tilde{V}}_{n}},  \qquad{{\tilde{C}}_{n}}{{x}^{p}}=\frac{\left( n+p+1 \right)!}{\left( p+1 \right)!}{{x}^{p}}.$$
{\bfseries Proof.} We use  Theorem 3 and the identities
$$\tilde{U}_{n}^{-1}\tilde{V}_{n}^{-1}{{x}^{p}}=\frac{n!}{\left( p+1 \right)!}{{\left( x-1 \right)}_{p}},  \qquad\tilde{U}_{n+p+1}^{-1}{{x}^{p}}={{\left( x-1 \right)}_{p}}{{\left[ x+1 \right]}_{n}}.$$
Then
$${{\tilde{F}}_{n}}\tilde{U}_{n}^{-1}\tilde{V}_{n}^{-1}{{x}^{p}}=\frac{\left( 2n \right)!}{n!}{{\tilde{U}}_{2n}}\frac{n!}{\left( p+1 \right)!}{{\left( x-1 \right)}_{p}}{{\left[ x+1 \right]}_{n}}=$$
$$=\frac{\left( 2n \right)!}{\left( p+1 \right)!}{{\left( 1-x \right)}^{n-p-1}}\frac{\left( n+p+1 \right)!}{\left( 2n \right)!}{{\tilde{U}}_{n+p+1}}{{\left( x-1 \right)}_{p}}{{\left[ x+1 \right]}_{n}}=$$
$$=\frac{\left( n+p+1 \right)!}{\left( p+1 \right)!}{{\left( 1-x \right)}^{n-p-1}}{{x}^{p}},$$
or
$${{\tilde{S}}_{n}}\tilde{V}_{n}^{-1}=\tilde{V}_{n}^{-1}{{\tilde{C}}_{n}}.$$

Thus,
$${{\tilde{S}}_{2}}=\left( \begin{matrix}
   1 & 0  \\
   -1 & 1  \\
\end{matrix} \right)\left( \begin{matrix}
   3! & 0  \\
   0 & {4!}/{2!}\;  \\
\end{matrix} \right)\left( \begin{matrix}
   1 & 0  \\
   1 & 1  \\
\end{matrix} \right),$$
$${{\tilde{S}}_{3}}=\left( \begin{matrix}
   1 & 0 & 0  \\
   -2 & 1 & 0  \\
   1 & -1 & 1  \\
\end{matrix} \right)\left( \begin{matrix}
   4! & 0 & 0  \\
   0 & {5!}/{2!}\; & 0  \\
   0 & 0 & {6!}/{3!}\;  \\
\end{matrix} \right)\left( \begin{matrix}
   1 & 0 & 0  \\
   2 & 1 & 0  \\
   1 & 1 & 1  \\
\end{matrix} \right),$$
$${{\tilde{S}}_{4}}=\left( \begin{matrix}
   1 & 0 & 0 & 0  \\
   -3 & 1 & 0 & 0  \\
   3 & -2 & 1 & 0  \\
   -1 & 1 & -1 & 1  \\
\end{matrix} \right)\left( \begin{matrix}
   5! & 0 & 0 & 0  \\
   0 & {6!}/{2!}\; & 0 & 0  \\
   0 & 0 & {7!}/{3!}\; & 0  \\
   0 & 0 & 0 & {8!}/{4!}\;  \\
\end{matrix} \right)\left( \begin{matrix}
   1 & 0 & 0 & 0  \\
   3 & 1 & 0 & 0  \\
   3 & 2 & 1 & 0  \\
   1 & 1 & 1 & 1  \\
\end{matrix} \right).$$

Note that since ${{\alpha }_{n}}\left( x \right)|{{\left( 1-x \right)}^{-1}}=x$, then
$$x{{\tilde{S}}_{n}}{{x}^{0}}={{\varphi }_{n}}\left( x \right)|{{\left( 1-x \right)}^{-1}}=\left( n+1 \right)!{{N}_{n}}\left( x \right);$$
since
$$x\frac{\left( 2n \right)!}{n!}\tilde{F}_{n}^{-1}{{x}^{0}}=x\left( x+n+1 \right)...\left( x+n+n-1 \right)=\left[ n,\to  \right]{{\left( 1,\log C\left( x \right) \right)}_{{{e}^{x}}}}$$
where $C\left( x \right)$ is the Catalan series, then
$$x\frac{\left( 2n \right)!}{n!}\tilde{S}_{n}^{-1}{{x}^{0}}={{\alpha }_{n}}\left( x \right)|C\left( x \right)=\frac{1}{n}\sum\limits_{m=1}^{n}{\left( \begin{matrix}
   -n  \\
   m-1  \\
\end{matrix} \right)\left( \begin{matrix}
   2n  \\
   n-m  \\
\end{matrix} \right){{x}^{m}}}.$$

Let $_{\left( \beta  \right)}a\left( x \right)$ is the generalized binomial series. Denote
$$[n,\searrow ]{{\left( 1,{{x}_{\left( \beta  \right)}}a\left( x \right) \right)}_{{{e}^{x}}}}=\frac{_{\left( \beta  \right)}{{\varphi }_{n}}\left( x \right)}{{{\left( 1-x \right)}^{2n+1}}}=\sum\limits_{m=0}^{\infty }{\frac{m}{m+\beta n}}\left( \begin{matrix}
   m+\beta n  \\
   n  \\
\end{matrix} \right){{\left[ m+1 \right]}_{n}}{{x}^{m}}.$$
{\bfseries Theorem 8.}
$$_{\left( \beta  \right)}{{\varphi }_{n}}\left( x \right)=\frac{\left( n+1 \right)!}{n}\sum\limits_{m=1}^{n}{\left( \begin{matrix}
   n\left( 2-\beta  \right)  \\
   m-1  \\
\end{matrix} \right)\left( \begin{matrix}
   n\beta   \\
   n-m  \\
\end{matrix} \right){{x}^{m}}}.$$
{\bfseries Proof.} Taking into account the polynomial argument, we prove the theorem for positive integers $\beta $. Since
$${{\tilde{V}}_{n}}{}_{\left( \beta  \right)}{{\tilde{\alpha }}_{n}}\left( x \right)={}_{\left( \beta  \right)}{{\tilde{v}}_{n}}\left( x \right)=\sum\limits_{m=0}^{n-1}{\frac{m+1}{n}\left( \begin{matrix}
   n\beta   \\
   n-m-1  \\
\end{matrix} \right)}{{x}^{m}},$$
$$\left[ m,\to  \right]\tilde{V}_{n}^{-1}=\sum\limits_{i=0}^{m}{\left( \begin{matrix}
   m-n  \\
   m-i  \\
\end{matrix} \right){{x}^{i}}}=\sum\limits_{i=0}^{m}{{{\left( -1 \right)}^{m-i}}}\left( \begin{matrix}
   n-i-1  \\
   m-i  \\
\end{matrix} \right){{x}^{i}},$$
then
$$\left[ {{x}^{m}} \right]{}_{\left( \beta  \right)}{{\tilde{\varphi }}_{n}}\left( x \right)=\left[ {{x}^{m}} \right]\tilde{V}_{n}^{-1}{{\tilde{C}}_{n}}{}_{\left( \beta  \right)}{{\tilde{v}}_{n}}\left( x \right)=$$
$$=\sum\limits_{i=0}^{m}{{{\left( -1 \right)}^{m-i}}}\left( \begin{matrix}
   n-i-1  \\
   m-i  \\
\end{matrix} \right)\frac{\left( n+i+1 \right)!}{\left( i+1 \right)!}\frac{\left( i+1 \right)}{n}\left( \begin{matrix}
   n\beta   \\
   n-i-1  \\
\end{matrix} \right)\frac{\left( n\beta -n+m+1 \right)!}{\left( n\beta -n+m+1 \right)!}=$$
$$=\frac{\left( n+1 \right)!}{n}\left( \begin{matrix}
   n\beta   \\
   n-m-1  \\
\end{matrix} \right)\sum\limits_{i=0}^{m}{{{\left( -1 \right)}^{m-i}}\left( \begin{matrix}
   n+i+1  \\
   i  \\
\end{matrix} \right)}\left( \begin{matrix}
   n\beta -n+m+1  \\
   m-i  \\
\end{matrix} \right)=$$
$$=\frac{\left( n+1 \right)!}{n}\left( \begin{matrix}
   n\beta   \\
   n-m-1  \\
\end{matrix} \right){{\left( -1 \right)}^{m}}\left( \begin{matrix}
   n\beta -2n+m-1  \\
   m  \\
\end{matrix} \right)=$$
$$=\frac{\left( n+1 \right)!}{n}\left( \begin{matrix}
   n\beta   \\
   n-m-1  \\
\end{matrix} \right)\left( \begin{matrix}
   n\left( 2-\beta  \right)  \\
   m  \\
\end{matrix} \right).$$

Note that
$$_{\left( 0 \right)}{{\varphi }_{n}}\left( x \right)=\frac{\left( 2n \right)!}{n!}{{x}^{n}},  \qquad_{\left( 1 \right)}{{\varphi }_{n}}\left( x \right)=\left( n+1 \right)!{{N}_{n}}\left( x \right),   \qquad_{\left( 2 \right)}{{\varphi }_{n}}\left( x \right)=\frac{\left( 2n \right)!}{n!}x;$$
since
$${{\left( 1,x{}_{\left( \beta  \right)}a\left( x \right) \right)}^{-1}}=\left( 1,{{x}_{\left( \beta -1 \right)}}{{a}^{-1}}\left( x \right) \right),  \qquad{}_{\left( \beta -1 \right)}{{a}^{-1}}\left( -x \right){{=}_{\left( 2-\beta  \right)}}a\left( x \right),$$ 
then
  $$_{\left( 2-\beta  \right)}{{\varphi }_{n}}\left( x \right)=x{{J}_{n}}{}_{\left( \beta  \right)}{{\varphi }_{n}}\left( x \right).$$
\section{Transformations of general form}

We introduce the matrices  ${{U}_{n}}$:
$${{U}_{n}}{{x}^{p}}=\frac{1}{n!}{{\left( 1-x \right)}^{n-p}}{{A}_{p}}\left( x \right),  \qquad{{A}_{0}}\left( x \right)=1,  \qquad p=0, \text{ }1,\text{ } … , \text{ } n.$$
Or
$${{U}_{n}}{{x}^{p}}={{\left( 1-x \right)}^{n+1}}\frac{1}{n!}\sum\limits_{m=0}^{\infty }{{{m}^{p}}}{{x}^{m}}.$$
For example,
$${{U}_{0}}=\left( 1 \right), \quad{{U}_{1}}=\left( \begin{matrix}
   1 & 0  \\
   -1 & 1  \\
\end{matrix} \right), \quad{{U}_{2}}=\frac{1}{2!}\left( \begin{matrix}
   1 & 0 & 0  \\
   -2 & 1 & 1  \\
   1 & -1 & 1  \\
\end{matrix} \right), \quad{{U}_{3}}=\frac{1}{3!}\left( \begin{matrix}
   1 & 0 & 0 & 0  \\
   -3 & 1 & 1 & 1  \\
   3 & -2 & 0 & 4  \\
   -1 & 1 & -1 & 1  \\
\end{matrix} \right).$$
Let
$$\left[ n,\to  \right]\left( b\left( x \right),\log a\left( x \right) \right)={{c}_{n}}\left( x \right)=\sum\limits_{m=0}^{n}{{{c}_{m}}{{x}^{m}}}, \qquad{{b}_{0}}\ne 0,$$
$$\left[ n,\to  \right]{{\left( b\left( x \right),\log a\left( x \right) \right)}_{{{e}^{x}}}}={{s}_{n}}\left( x \right)=\sum\limits_{m=0}^{n}{{{s}_{m}}{{x}^{m}}},$$
$$\left[ n,\to  \right]\left( b\left( x \right),a\left( x \right) \right)=\frac{{{g}_{n}}\left( x \right)}{{{\left( 1-x \right)}^{n+1}}}.$$
Since
$$\left( b\left( x \right),a\left( x \right) \right)=\left( b\left( x \right),\log a\left( x \right) \right)\left( 1,{{e}^{x}} \right),$$
then
$$\frac{{{g}_{n}}\left( x \right)}{{{\left( 1-x \right)}^{n+1}}}=\sum\limits_{p=0}^{n}{\frac{{{{c}_{p}}{{A}_{p}}\left( x \right)}/{p!}\;}{{{\left( 1-x \right)}^{p+1}}}}=\frac{1}{n!}\sum\limits_{p=0}^{n}{\frac{{{s}_{p}}{{A}_{p}}\left( x \right)}{{{\left( 1-x \right)}^{p+1}}}}=\frac{\frac{1}{n!}\sum\limits_{p=0}^{n}{{{s}_{p}}{{\left( 1-x \right)}^{n-p}}{{A}_{p}}\left( x \right)}}{{{\left( 1-x \right)}^{n+1}}}.$$
Thus,
$$b\left( x \right){{a}^{m}}\left( x \right)=\sum\limits_{n=0}^{\infty }{\frac{{{s}_{n}}\left( m \right)}{n!}{{x}^{n}}},  \qquad\frac{{{g}_{n}}\left( x \right)}{{{\left( 1-x \right)}^{n+1}}}=\sum\limits_{m=0}^{\infty }{\frac{{{s}_{n}}\left( m \right)}{n!}{{x}^{m}}},  \qquad{{g}_{n}}\left( x \right)={{U}_{n}}{{s}_{n}}\left( x \right).$$
Since
$$\frac{1}{{{\left( 1-x \right)}^{n+1}}}=\sum\limits_{m=0}^{\infty }{\frac{{{\left[ m+1 \right]}_{n}}}{n!}}{{x}^{m}},$$
then
$$U_{n}^{-1}{{x}^{0}}={{\left[ x+1 \right]}_{n}},  \qquad U_{n}^{-1}{{x}^{p}}=x\tilde{U}_{n}^{-1}{{x}^{p-1}},$$
or
$$U_{n}^{-1}{{x}^{p}}={{\left( x \right)}_{p}}{{\left[ x+1 \right]}_{n-p}}.$$
For example,
$$U_{3}^{-1}=\left( \begin{matrix}
   6 & 0 & 0 & 0  \\
   11 & 2 & -1 & 2  \\
   6 & 3 & 0 & -3  \\
   1 & 1 & 1 & 1  \\
\end{matrix} \right).$$
{\bfseries Theorem 9.}
$${{U}_{n}}E\left( 1,-x \right)U_{n}^{-1}={{\left( -1 \right)}^{n}}{{J}_{n}}.$$
{\bfseries Proof.}
$$E\left( 1,-x \right){{\left( x \right)}_{p}}{{\left[ x+1 \right]}_{n-p}}={{\left( -x-1 \right)}_{p}}{{\left[ -x \right]}_{n-p}}={{\left( -1 \right)}^{n}}{{\left( x \right)}_{n-p}}{{\left[ x+1 \right]}_{p}},$$
or
$$E\left( 1,-x \right)U_{n}^{-1}={{\left( -1 \right)}^{n}}U_{n}^{-1}{{J}_{n}}.$$

Thus,
$${{\left( -1 \right)}^{n}}{{J}_{n}}{{g}_{n}}\left( x \right)={{U}_{n}}{{s}_{n}}\left( -x-1 \right).$$
Since
$$\left( b\left( x \right),\log a\left( x \right) \right)\left( 1,-x \right)\left( {{e}^{x}},x \right)=\left( b\left( x \right){{a}^{-1}}\left( x \right),\log {{a}^{-1}}\left( x \right) \right),$$
then polynomials  ${{\left( -1 \right)}^{n}}{{J}_{n}}{{g}_{n}}\left( x \right)$ are the numerator polynomials of the matrix 
$$\left( b\left( x \right){{a}^{-1}}\left( x \right),x{{a}^{-1}}\left( x \right) \right).$$
In particular,
$${{\left( -1 \right)}^{n}}{{J}_{n}}{{\alpha }_{n}}\left( x \right)=\tilde{\alpha }_{n}^{\left( -1 \right)}\left( x \right).$$

Denote
$${{V}_{n}}={{J}_{n}}E{{J}_{n}}=\left( {{\left( 1+x \right)}^{n+1}},x \right){{P}^{-1}}{{I}_{n}}={{\tilde{V}}_{n+1}},$$
$$\left[ n,\to  \right]\left( b\left( x \right),a\left( x \right)-1 \right)={{w}_{n}}\left( x \right).$$
Since
$$\left( b\left( x \right),a\left( x \right) \right)=\left( b\left( x \right),a\left( x \right)-1 \right)\left( 1,1+x \right),$$
then
$${{g}_{n}}\left( x \right)=V_{n}^{-1}{{w}_{n}}\left( x \right).$$
{\bfseries Theorem 10.}\emph{ If ${{s}_{n-m}}\left( x \right)$ is a polynomial of degree $n-m$, then
$${{U}_{n}}{{s}_{n-m}}\left( x \right)={{\left( 1-x \right)}^{m}}\frac{\left( n-m \right)!}{n!}{{U}_{n-m}}{{s}_{n-m}}\left( x \right).$$
Respectively, if  ${{c}_{n-m}}\left( x \right)$ is a polynomial of degree $\le n-m$, then
$$U_{n}^{-1}{{\left( 1-x \right)}^{m}}{{c}_{n-m}}\left( x \right)=\frac{n!}{\left( n-m \right)!}U_{n-m}^{-1}{{c}_{n-m}}\left( x \right).$$}
{\bfseries Proof.}
$$\left( {{\left( 1-x \right)}^{-m}},x \right){{U}_{n}}{{I}_{n-m}}=\frac{\left( n-m \right)!}{n!}{{U}_{n-m}},$$
$$U_{n}^{-1}\left( {{\left( 1-x \right)}^{m}},x \right){{I}_{n-m}}=\frac{n!}{\left( n-m \right)!}U_{n-m}^{-1}.$$

From this we find:
$$U_{n}^{-1}V_{n}^{-1}{{x}^{p}}=U_{n}^{-1}{{\left( 1-x \right)}^{n-p}}{{x}^{p}}=\frac{n!}{p!}{{\left( x \right)}_{p}}=\frac{n!}{p!}\sum\limits_{m=0}^{p}{s\left( p,\text{ }m \right){{x}^{m}}},$$
$${{V}_{n}}{{U}_{n}}{{x}^{p}}=\frac{1}{n!}\sum\limits_{m=0}^{p}{m!S\left( p,\text{ }m \right)}\text{ }{{x}^{m}}.$$
{\bfseries Remark 1.} Matrices ${{U}_{n}}$, $U_{n}^{-1}$ are  associated with the matrices ${{\tilde{U}}_{n}}$, $\tilde{U}_{n}^{-1}$ by the identities
$${{\tilde{U}}_{n}}={{\left( x,x \right)}^{T}}{{U}_{n}}\left( x,x \right),   \qquad\tilde{U}_{n}^{-1}={{\left( x,x \right)}^{T}}U_{n}^{-1}\left( x,x \right),$$
and, since
$$\frac{{{{\tilde{\alpha }}}_{n}}\left( x \right)}{{{\left( 1-x \right)}^{n+1}}}=\sum\limits_{m=0}^{\infty }{\frac{{{u}_{n}}\left( m+1 \right)}{n!}}{{x}^{m}},  \qquad{{\tilde{\alpha }}_{n}}\left( x \right)={{U}_{n}}E{{u}_{n}}\left( x \right),$$ 
then
 $${{\tilde{U}}_{n}}={{U}_{n}}E\left( x,x \right){{I}_{n-1}},   \qquad\tilde{U}_{n}^{-1}={{\left( x,x \right)}^{T}}{{E}^{-1}}U_{n}^{-1}{{I}_{n-1}}.$$

Denote
$$[n,\searrow ]{{\left( b\left( x \right),xa\left( x \right) \right)}_{{{e}^{x}}}}=\sum\limits_{m=0}^{\infty }{\frac{{{\left[ m+1 \right]}_{n}}{{s}_{n}}\left( m \right)}{n!}{{x}^{m}}}=\frac{{{h}_{n}}\left( x \right)}{{{\left( 1-x \right)}^{2n+1}}},$$
Then
$${{h}_{n}}\left( x \right)=\frac{\left( 2n \right)!}{n!}{{U}_{2n}}{{\left[ x+1 \right]}_{n}}{{s}_{n}}\left( x \right).$$
Polynomials ${{g}_{n}}\left( x \right)$, ${{h}_{n}}\left( x \right)$ will be called respectively the ordinary and the exponential numerator polynomials. Names GEP and GNP we will fix for the polynomials ${{\alpha }_{n}}\left( x \right)$, ${{\varphi }_{n}}\left( x \right)$.

We introduce the matrices ${{F}_{n}}$:
$${{F}_{n}}=\frac{\left( 2n \right)!}{n!}{{U}_{2n}}\left( {{\left[ x+1 \right]}_{n}},x \right){{I}_{n}},  \qquad F_{n}^{-1}=\frac{n!}{\left( 2n \right)!}{{\left( {{\left[ x+1 \right]}_{n}},x \right)}^{-1}}U_{2n}^{-1}{{I}_{n}};$$
Since
$${{U}_{2n}}{{\left[ x+1 \right]}_{n}}={{\left( 1-x \right)}^{n}}\frac{n!}{\left( 2n \right)!}{{U}_{n}}{{\left[ x+1 \right]}_{n}}=\frac{n!}{\left( 2n \right)!}{{\left( 1-x \right)}^{n}},$$
then
$${{F}_{n}}{{x}^{0}}={{\left( 1-x \right)}^{n}}, \qquad{{F}_{n}}{{x}^{p}}=x{{\tilde{F}}_{n}}{{x}^{p-1}},$$
or
$${{F}_{n}}{{x}^{p}}={{\left( 1-x \right)}^{2n+1}}\sum\limits_{m=0}^{\infty }{{{m}^{p}}}\left( \begin{matrix}
   m+n  \\
   n  \\
\end{matrix} \right){{x}^{m}}.$$
For example,
$${{F}_{2}}=\left( \begin{matrix}
   1 & 0 & 0  \\
   -2 & 3 & 3  \\
   1 & -3 & 9  \\
\end{matrix} \right), \qquad{{F}_{3}}=\left( \begin{matrix}
   1 & 0 & 0 & 0  \\
   -3 & 4 & 4 & 4  \\
   3 & -8 & 12 & 52  \\
   -1 & 4 & -16 & 64  \\
\end{matrix} \right).$$
Respectively,
$$F_{n}^{-1}{{x}^{0}}=\frac{n!}{\left( 2n \right)!}{{\left[ x+n+1 \right]}_{n}}, \qquad F_{n}^{-1}{{x}^{p}}=x\tilde{F}_{n}^{-1}{{x}^{p-1}},$$
or
$$F_{n}^{-1}{{x}^{p}}=\frac{n!}{\left( 2n \right)!}{{\left( x \right)}_{p}}{{\left[ x+n+1 \right]}_{n-p}}.$$
For example,
$$F_{2}^{-1}=\frac{2!}{4!}\left( \begin{matrix}
   12 & 0 & 0  \\
   7 & 3 & -1  \\
   1 & 1 & 1  \\
\end{matrix} \right), \qquad F_{3}^{-1}=\frac{3!}{6!}\left( \begin{matrix}
   120 & 0 & 0 & 0  \\
   74 & 20 & -4 & 2  \\
   15 & 9 & 3 & -3  \\
   1 & 1 & 1 & 1  \\
\end{matrix} \right).$$
Thus,
$${{F}_{n}}{{s}_{n}}\left( x \right)={{h}_{n}}\left( x \right).$$
{\bfseries Example 1.} We explain the identity
$${{F}_{n}}{{\left[ x+n+1 \right]}_{n}}=\frac{\left( 2n \right)!}{n!}.$$
Analog of the identity for ordinary Riordan arrays 
$$\left[ n,\searrow  \right]\left( a\left( x \right),xa\left( x \right) \right)=\frac{{{{\tilde{\alpha }}}_{n}}\left( x \right)}{{{\left( 1-x \right)}^{n+1}}}$$ 
is the identity for exponential Riordan arrays
$$\left[ n,\searrow  \right]{{\left( {{\left( xa\left( x \right) \right)}^{\prime }},xa\left( x \right) \right)}_{{{e}^{x}}}}=\frac{{{{\tilde{\varphi }}}_{n}}\left( x \right)}{{{\left( 1-x \right)}^{2n+1}}}.$$
Since
$$\left[ n,\to  \right]{{\left( 1+x{{\left( \log a\left( x \right) \right)}^{\prime }},\log a\left( x \right) \right)}_{{{e}^{x}}}}=\left( x+n \right){{\tilde{u}}_{n}}\left( x \right),$$
$$\left( {{\left( xa\left( x \right) \right)}^{\prime }},\log a\left( x \right) \right)=\left( 1+x{{\left( \log a\left( x \right) \right)}^{\prime }},\log a\left( x \right) \right)\left( {{e}^{x}},x \right),$$
then
$$\left[ n,\to  \right]{{\left( {{\left( xa\left( x \right) \right)}^{\prime }},\log a\left( x \right) \right)}_{{{e}^{x}}}}=\left( x+n+1 \right){{\tilde{u}}_{n}}\left( x+1 \right).$$
If $a\left( x \right)=C\left( x \right)$, then
$${{\tilde{\varphi }}_{n}}\left( x \right)=\frac{\left( 2n \right)!}{n!},  \qquad{{\tilde{u}}_{n}}\left( x \right)={{\left[ x+n+1 \right]}_{n-1}},  \qquad\left( x+n+1 \right){{\tilde{u}}_{n}}\left( x+1 \right)={{\left[ x+n+1 \right]}_{n}}.$$
{\bfseries Theorem 11.}
$${{F}_{n}}{{E}^{n+1}}\left( 1,-x \right)F_{n}^{-1}={{\left( -1 \right)}^{n}}{{J}_{n}}.$$
{\bfseries Proof.}
$${{E}^{n+1}}\left( 1,-x \right){{\left( x \right)}_{p}}{{\left[ x+n+1 \right]}_{n-p}}={{\left( -x-n-1 \right)}_{p}}{{\left[ -x \right]}_{n-p}}={{\left( -1 \right)}^{n}}{{\left( x \right)}_{n-p}}{{\left[ x+n+1 \right]}_{p}},$$
or
$${{E}^{n+1}}\left( 1,-x \right)F_{n}^{-1}={{\left( -1 \right)}^{n}}F_{n}^{-1}{{J}_{n}}.$$

Thus,
$${{\left( -1 \right)}^{n}}{{J}_{n}}{{h}_{n}}\left( x \right)={{F}_{n}}E{{s}_{n}}\left( -x-n \right).$$
Since (see Remark 2)
$${{s}_{n}}\left( -x-n \right)=\left[ n,\to  \right]{{\left( b\left( x{}_{\left( -1 \right)}{{a}^{-1}}\left( x \right) \right)\left( 1+x{{\left( \log {}_{\left( -1 \right)}{{a}^{-1}}\left( x \right) \right)}^{\prime }} \right),\log {}_{\left( -1 \right)}{{a}^{-1}}\left( x \right) \right)}_{{{e}^{x}}}},$$
where
$$\left( 1,x{}_{\left( -1 \right)}{{a}^{-1}}\left( x \right) \right)={{\left( 1,xa\left( x \right) \right)}^{-1}},  \qquad_{\left( -1 \right)}{{a}^{\varphi }}\left( x \right)=\sum\limits_{m=0}^{\infty }{\frac{\varphi }{\left( \varphi -n \right)}\frac{{{u}_{n}}\left( \varphi -n \right)}{n!}{{x}^{m}}},$$
then polynomials ${{\left( -1 \right)}^{n}}{{J}_{n}}{{h}_{n}}\left( x \right)$ are the numerator polynomials of the matrix 
$${{\left( b\left( x{}_{\left( -1 \right)}{{a}^{-1}}\left( x \right) \right){{\left( x{}_{\left( -1 \right)}{{a}^{-1}}\left( x \right) \right)}^{\prime }},x{}_{\left( -1 \right)}{{a}^{-1}}\left( x \right) \right)}_{{{e}^{x}}}}.$$
In particular,
$${{\left( -1 \right)}^{n}}{{J}_{n}}{{\varphi }_{n}}\left( x \right)=\tilde{\varphi }_{n}^{\left[ -1 \right]}\left( x \right).$$
{\bfseries Remark 2.} We represent the matrix ${{\left( b\left( x \right),{{a}^{-1}}\left( x \right) \right)}^{T}}$ in the form
$${{\left( b\left( x \right),{{a}^{-1}}\left( x \right) \right)}^{T}}=\left( \begin{matrix}
   {{s}_{0}}\left( 0 \right) & {{s}_{1}}\left( 0 \right) & {{s}_{2}}\left( 0 \right) & {{s}_{3}}\left( 0 \right) & \cdots   \\
   {{s}_{0}}\left( -1 \right) & {{s}_{1}}\left( -1 \right) & {{s}_{2}}\left( -1 \right) & {{s}_{3}}\left( -1 \right) & \cdots   \\
   {{s}_{0}}\left( -2 \right) & {{s}_{1}}\left( -2 \right) & {{s}_{2}}\left( -2 \right) & {{s}_{3}}\left( -2 \right) & \cdots   \\
   {{s}_{0}}\left( -3 \right) & {{s}_{1}}\left( -3 \right) & {{s}_{2}}\left( -3 \right) & {{s}_{3}}\left( -3 \right) & \cdots   \\
   \vdots  & \vdots  & \vdots  & \vdots  & \ddots   \\
\end{matrix} \right)\left| {{e}^{x}} \right|,$$
where ${{s}_{0}}\left( x \right)={{b}_{0}}$. From the Lagrange inversion theorem it follows that
$$\left[ n,\searrow  \right]{{\left( b\left( x \right),{{a}^{-1}}\left( x \right) \right)}^{T}}=b\left( x{}_{\left( -1 \right)}{{a}^{-1}}\left( x \right) \right)\left( 1+x{{\left( \log {}_{\left( -1 \right)}{{a}^{-1}}\left( x \right) \right)}^{\prime }} \right){}_{\left( -1 \right)}{{a}^{-n}}\left( x \right).$$

We introduce the matrices ${{S}_{n}}$:
$${{S}_{n}}={{F}_{n}}U_{n}^{-1},  \qquad S_{n}^{-1}={{U}_{n}}F_{n}^{-1}.$$
For example,
$${{S}_{2}}=2!\left( \begin{matrix}
   1 & 0 & 0  \\
   4 & 3 & 0  \\
   1 & 3 & 6  \\
\end{matrix} \right), \quad{{S}_{3}}=3!\left( \begin{matrix}
   1 & 0 & 0 & 0  \\
   9 & 4 & 0 & 0  \\
   9 & 12 & 10 & 0  \\
   1 & 4 & 10 & 20  \\
\end{matrix} \right),
\quad{{S}_{4}}=4!\left( \begin{matrix}
   1 & 0 & 0 & 0 & 0  \\
   16 & 5 & 0 & 0 & 0  \\
   36 & 30 & 15 & 0 & 0  \\
   16 & 30 & 40 & 35 & 0  \\
   1 & 5 & 15 & 35 & 70  \\
\end{matrix} \right);$$
$$S_{2}^{-1}=\frac{2!}{4!}\left( \begin{matrix}
   6 & 0 & 0  \\
   -8 & 2 & 0  \\
   3 & -1 & 1  \\
\end{matrix} \right), \quad S_{3}^{-1}=\frac{3!}{6!}\left( \begin{matrix}
   20 & 0 & 0 & 0  \\
   -45 & 5 & 0 & 0  \\
   36 & -6 & 2 & 0  \\
   -10 & 2 & -1 & 1  \\
\end{matrix} \right),$$
$$ S_{4}^{-1}=\frac{4!}{8!}\left( \begin{matrix}
   70 & 0 & 0 & 0 & 0  \\
   -224 & 14 & 0 & 0 & 0  \\
   280 & -28 & {14}/{3}\; & 0 & 0  \\
   -160 & 20 & {-16}/{3}\; & 2 & 0  \\
   35 & -5 & {5}/{3}\; & -1 & 1  \\
\end{matrix} \right).$$
Then
$${{S}_{n}}{{g}_{n}}\left( x \right)={{h}_{n}}\left( x \right).$$
{\bfseries Theorem 12.}
$${{S}_{n}}=V_{n}^{-1}{{C}_{n}}{{V}_{n}},  \qquad{{C}_{n}}{{x}^{p}}=\frac{\left( n+p \right)!}{p!}{{x}^{p}}.$$
{\bfseries Proof.} We use  Theorem 10 and the identities
$$U_{n}^{-1}V_{n}^{-1}{{x}^{p}}=\frac{n!}{p!}{{\left( x \right)}_{p}},  \qquad U_{n+p}^{-1}{{x}^{p}}={{\left( x \right)}_{p}}{{\left[ x+1 \right]}_{n}}.$$
Then
$${{F}_{n}}U_{n}^{-1}V_{n}^{-1}{{x}^{p}}=\frac{\left( 2n \right)!}{n!}{{U}_{2n}}\frac{n!}{p!}{{\left( x \right)}_{p}}{{\left[ x+1 \right]}_{n}}=$$
$$=\frac{\left( 2n \right)!}{p!}{{\left( 1-x \right)}^{n-p}}\frac{\left( n+p \right)!}{\left( 2n \right)!}{{U}_{n+p}}{{\left( x \right)}_{p}}{{\left[ x+1 \right]}_{n}}=\frac{\left( n+p \right)!}{p!}{{\left( 1-x \right)}^{n-p}}{{x}^{p}},$$
or
$${{S}_{n}}V_{n}^{-1}=V_{n}^{-1}{{C}_{n}}.$$
{\bfseries Theorem 13.}
$${{S}_{n}}{{x}^{p}}=\frac{\left( n+p \right)!\left( n-p \right)!}{n!}\sum\limits_{m=p}^{n}{\left( \begin{matrix}
   n  \\
   m-p  \\
\end{matrix} \right)\left( \begin{matrix}
   n  \\
   n-m  \\
\end{matrix} \right){{x}^{m}}}.$$
{\bfseries Proof.}
$$\left[ m,\to  \right]V_{n}^{-1}=\sum\limits_{i=0}^{m}{\left( \begin{matrix}
   m-n-1  \\
   m-i  \\
\end{matrix} \right){{x}^{i}}}=\sum\limits_{i=0}^{m}{{{\left( -1 \right)}^{m-i}}}\left( \begin{matrix}
   n-i  \\
   m-i  \\
\end{matrix} \right){{x}^{i}},$$
$${{C}_{n}}{{V}_{n}}{{x}^{p}}=\sum\limits_{i=p}^{n}{\frac{\left( n+i \right)!}{i!}\left( \begin{matrix}
   n-p  \\
   i-p  \\
\end{matrix} \right)}{{x}^{i}},$$
$$\left[ {{x}^{m}} \right]V_{n}^{-1}{{C}_{n}}{{V}_{n}}{{x}^{p}}=\sum\limits_{i=p}^{m}{{{\left( -1 \right)}^{m-i}}\left( \begin{matrix}
   n-i  \\
   m-i  \\
\end{matrix} \right)\frac{\left( n+i \right)!}{i!}\left( \begin{matrix}
   n-p  \\
   i-p  \\
\end{matrix} \right)}=$$
$$=\frac{\left( n+p \right)!\left( n-p \right)!}{\left( n-m \right)!m!}\sum\limits_{i=p}^{m}{{{\left( -1 \right)}^{m-i}}}\left( \begin{matrix}
   n+i  \\
   i-p  \\
\end{matrix} \right)\left( \begin{matrix}
   m  \\
   i  \\
\end{matrix} \right)=$$ 
$$=\frac{\left( n+p \right)!\left( n-p \right)!}{\left( n-m \right)!m!}{{\left( -1 \right)}^{m-p}}\left( \begin{matrix}
   m-n-p-1  \\
   m-p  \\
\end{matrix} \right)=$$
$$=\frac{\left( n+p \right)!\left( n-p \right)!}{n!}\left( \begin{matrix}
   n  \\
   m-p  \\
\end{matrix} \right)\left( \begin{matrix}
   n  \\
   n-m  \\
\end{matrix} \right).$$
{\bfseries Theorem 14.}
$$S_{n}^{-1}{{x}^{p}}=\frac{p!\left( n-p \right)!}{\left( 2n \right)!}\sum\limits_{m=p}^{n}{\left( \begin{matrix}
   -n  \\
   m-p  \\
\end{matrix} \right)\left( \begin{matrix}
   2n  \\
   n-m  \\
\end{matrix} \right){{x}^{m}}}.$$
{\bfseries Proof.}
$$\left[ {{x}^{m}} \right]V_{n}^{-1}C_{n}^{-1}{{V}_{n}}{{x}^{p}}=\sum\limits_{i=p}^{m}{{{\left( -1 \right)}^{m-i}}\left( \begin{matrix}
   n-i  \\
   m-i  \\
\end{matrix} \right)\frac{i!}{\left( n+i \right)!}\left( \begin{matrix}
   n-p  \\
   i-p  \\
\end{matrix} \right)}=$$
$$=\frac{p!\left( n-p \right)!}{\left( n-m \right)!\left( n+m \right)!}\sum\limits_{i=p}^{m}{{{\left( -1 \right)}^{m-i}}}\left( \begin{matrix}
   i  \\
   i-p  \\
\end{matrix} \right)\left( \begin{matrix}
   n+m  \\
   m-i  \\
\end{matrix} \right)=$$ 
$$=\frac{p!\left( n-p \right)!}{\left( n-m \right)!\left( n+m \right)!}{{\left( -1 \right)}^{m-p}}\left( \begin{matrix}
   n+m-p-1  \\
   m-p  \\
\end{matrix} \right)=$$
$$=\frac{p!\left( n-p \right)!}{\left( 2n \right)!}\left( \begin{matrix}
   -n  \\
   m-p  \\
\end{matrix} \right)\left( \begin{matrix}
   2n  \\
   n-m  \\
\end{matrix} \right).$$
\section{Generalized Narayana polynomials of type B}
Polynomials 
$${}^{B}{{N}_{n}}\left( x \right)=\sum\limits_{m=0}^{n}{{{\left( \begin{matrix}
   n  \\
   m  \\
\end{matrix} \right)}^{2}}{{x}^{m}}}={{\left( 1-x \right)}^{2n+1}}\sum\limits_{m=0}^{\infty }{{{\left( \begin{matrix}
   m+n  \\
   n  \\
\end{matrix} \right)}^{2}}{{x}^{m}}}$$
are called Narayana polynomials of type B. Denote
$$[n,\searrow ]{{\left( a\left( x \right),xa\left( x \right) \right)}_{{{e}^{x}}}}=\frac{{}^{B}{{\varphi }_{n}}\left( x \right)}{{{\left( 1-x \right)}^{2n+1}}}.$$
Let ${{\tilde{\alpha }}_{n}}\left( x \right)|a\left( x \right)$, $^{B}{{\varphi }_{n}}\left( x \right)|a\left( x \right)$ denotes respectively polynomials ${{\tilde{\alpha }}_{n}}\left( x \right)$, $^{B}{{\varphi }_{n}}\left( x \right)$,  associated with the matrices $\left( a\left( x \right),xa\left( x \right) \right)$, ${{\left( a\left( x \right),xa\left( x \right) \right)}_{{{e}^{x}}}}$. Then
$$^{B}{{\varphi }_{n}}\left( x \right)|{{\left( 1-x \right)}^{-1}}=n!{}^{B}{{N}_{n}}\left( x \right).$$
In this connection we will called polynomials $^{B}{{\varphi }_{n}}\left( x \right)$ the generalized Narayana polynomials of type B. Since
$$[n,\searrow ]\left( a\left( x \right),xa\left( x \right) \right)=\frac{{{{\tilde{\alpha }}}_{n}}\left( x \right)}{{{\left( 1-x \right)}^{n+1}}}=\sum\limits_{m=0}^{\infty }{\frac{{{u}_{n}}\left( m+1 \right)}{n!}}{{x}^{m}},$$ 
then
$$[n,\searrow ]{{\left( a\left( x \right),xa\left( x \right) \right)}_{{{e}^{x}}}}=\sum\limits_{m=0}^{\infty }{\frac{{{\left[ m+1 \right]}_{n}}{{u}_{n}}\left( m+1 \right)}{n!}{{x}^{m}}},$$
$$^{B}{{\varphi }_{n}}\left( x \right)=\frac{\left( 2n \right)!}{n!}{{U}_{2n}}{{\left[ x+1 \right]}_{n}}{{u}_{n}}\left( x+1 \right)={{F}_{n}}{{u}_{n}}\left( x+1 \right).$$
We introduce the matrices ${}^{B}{{F}_{n}}={{F}_{n}}E$:
$${}^{B}{{F}_{n}}{{x}^{p}}={{\left( 1-x \right)}^{2n+1}}\sum\limits_{m=0}^{\infty }{{{\left( m+1 \right)}^{p}}}\left( \begin{matrix}
   m+n  \\
   n  \\
\end{matrix} \right){{x}^{m}},\quad
^{B}F_{n}^{-1}{{x}^{p}}=\frac{n!}{\left( 2n \right)!}{{\left( x-1 \right)}_{p}}{{\left[ x+n \right]}_{n-p}}.$$
For example,
$$^{B}{{F}_{1}}=\left( \begin{matrix}
   1 & 1  \\
   -1 & 1  \\
\end{matrix} \right), \qquad^{B}{{F}_{2}}=\left( \begin{matrix}
   1 & 1 & 1  \\
   -2 & 1 & 7  \\
   1 & -2 & 4  \\
\end{matrix} \right), \qquad^{B}{{F}_{3}}=\left( \begin{matrix}
   1 & 1 & 1 & 1  \\
   -3 & 1 & 9 & 25  \\
   3 & -5 & -1 & 67  \\
   -1 & 3 & -9 & 27  \\
\end{matrix} \right),$$
$$^{B}F_{1}^{-1}=\frac{1}{2}\left( \begin{matrix}
   1 & -1  \\
   1 & 1  \\
\end{matrix} \right), \quad^{B}F_{2}^{-1}=\frac{2!}{4!}\left( \begin{matrix}
   6 & -2 & 2  \\
   5 & 1 & -3  \\
   1 & 1 & 1  \\
\end{matrix} \right), \quad^{B}F_{3}^{-1}=\frac{3!}{6!}\left( \begin{matrix}
   60 & -12 & 6 & -6  \\
   47 & 5 & -7 & 11  \\
   12 & 6 & 0 & -6  \\
   1 & 1 & 1 & 1  \\
\end{matrix} \right).$$
Then
$$^{B}{{F}_{n}}{{u}_{n}}\left( x \right)={}^{B}{{\varphi }_{n}}\left( x \right).$$
In particular,
$$^{B}{{F}_{n}}{{x}^{n}}={}^{B}{{\varphi }_{n}}\left( x \right)|{{e}^{x}}.$$

For the matrices $^{B}{{F}_{n}}$, Theorem 11 takes the simpler form. Since
$${{E}^{n-1}}\left( 1,-x \right){{\left( x-1 \right)}_{p}}{{\left[ x+n \right]}_{n-p}}={{\left( -x-n \right)}_{p}}{{\left[ -x+1 \right]}_{n-p}}={{\left( -1 \right)}^{n}}{{\left( x-1 \right)}_{n-p}}{{\left[ x+n \right]}_{p}},$$
then
$${}^{B}{{F}_{n}}{{E}^{n-1}}\left( 1,-x \right){}^{B}F_{n}^{-1}={{\left( -1 \right)}^{n}}{{J}_{n}},   \qquad{{\left( -1 \right)}^{n}}{{J}_{n}}{}^{B}{{\varphi }_{n}}\left( x \right)={{F}_{n}}{{u}_{n}}\left( -x-n \right),$$
where
$${{u}_{n}}\left( -x-n \right)=\left[ n,\to  \right]{{\left( 1+x{{\left( \log {}_{\left( -1 \right)}{{a}^{-1}}\left( x \right) \right)}^{\prime }},\log {}_{\left( -1 \right)}{{a}^{-1}}\left( x \right) \right)}_{{{e}^{x}}}}.$$
Denote
$$\left[ n,\searrow  \right]\left( 1+x{{\left( \log a\left( x \right) \right)}^{\prime }},xa\left( x \right) \right)_{{{e}^{x}}}^{-1}=\frac{^{B}\varphi _{n}^{\left[ -1 \right]}\left( x \right)}{{{\left( 1-x \right)}^{2n+1}}}.$$
Since 
$$\left( 1+x{{\left( \log a\left( x \right) \right)}^{\prime }},xa\left( x \right) \right)_{{{e}^{x}}}^{-1}={{\left( 1+x{{\left( \log {}_{\left( -1 \right)}{{a}^{-1}}\left( x \right) \right)}^{\prime }},x{}_{\left( -1 \right)}{{a}^{-1}}\left( x \right) \right)}_{{{e}^{x}}}},$$
then
$${}^{B}\varphi _{n}^{\left[ -1 \right]}\left( x \right)={{\left( -1 \right)}^{n}}{{J}_{n}}{}^{B}{{\varphi }_{n}}\left( x \right).$$

Let all polynomials $^{B}{{\varphi }_{n}}\left( x \right)$ are symmetric, i.e.
$${}^{B}{{\varphi }_{n}}\left( x \right)={{J}_{n}}{}^{B}{{\varphi }_{n}}\left( x \right).$$
Then ${{\left( xa\left( x \right) \right)}^{\prime }}={{a}^{2}}\left( x \right)$, or $\sum\nolimits_{m=0}^{n}{{{a}_{n-m}}{{a}_{m}}}=\left( n+1 \right){{a}_{n}}$. This is possible only in the case $a\left( x \right)={{\left( 1-\beta x \right)}^{-1}}$:
$$^{B}{{\varphi }_{n}}|{{\left( 1-\beta x \right)}^{-1}}={{\beta }^{n}}n!{}^{B}{{N}_{n}}\left( x \right).$$
{\bfseries Example 2.} Since 
$${{\tilde{\alpha }}_{n}}\left( x \right)|1+x={{x}^{n-1}},$$
then 
$$^{B}{{\varphi }_{n}}\left( x \right)|1+x={{S}_{n}}{{x}^{n-1}}=\frac{\left( 2n \right)!}{n!2}\left( 1+x \right){{x}^{n-1}};$$
$${{\left( 1+x,x\left( 1+x \right) \right)}_{{{e}^{x}}}}={{\left| {{e}^{x}} \right|}^{-1}}\left( \begin{matrix}
   1 & 0 & 0 & 0 & 0 & 0 & 0 & \cdots   \\
   1 & 1 & 0 & 0 & 0 & 0 & 0 & \cdots   \\
   0 & 2 & 1 & 0 & 0 & 0 & 0 & \cdots   \\
   0 & 1 & 3 & 1 & 0 & 0 & 0 & \cdots   \\
   0 & 0 & 3 & 4 & 1 & 0 & 0 & \cdots   \\
   0 & 0 & 1 & 6 & 5 & 1 & 0 & \cdots   \\
   0 & 0 & 0 & 4 & 10 & 6 & 1 & \cdots   \\
   \vdots  & \vdots  & \vdots  & \vdots  & \vdots  & \vdots  & \vdots  & \ddots   \\
\end{matrix} \right)\left| {{e}^{x}} \right|.$$
If $a\left( x \right)=1+x$, then $_{\left( -1 \right)}{{a}^{-1}}\left( x \right)=C\left( -x \right)$ and, hence,
$$\left[ n,\searrow  \right]{{\left( 1+x{{\left( \log C\left( x \right) \right)}^{\prime }},xC\left( x \right) \right)}_{{{e}^{x}}}}=\frac{\left( 2n \right)!}{n!2}\frac{1+x}{{{\left( 1-x \right)}^{2n+1}}};$$
$${{\left( 1+x{{\left( \log C\left( x \right) \right)}^{\prime }},xC\left( x \right) \right)}_{{{e}^{x}}}}={{\left| {{e}^{x}} \right|}^{-1}}\left( \begin{matrix}
   1 & 0 & 0 & 0 & 0 & 0 & \cdots   \\
   1 & 1 & 0 & 0 & 0 & 0 & \cdots   \\
   3 & 2 & 1 & 0 & 0 & 0 & \cdots   \\
   10 & 6 & 3 & 1 & 0 & 0 & \cdots   \\
   35 & 20 & 10 & 4 & 1 & 0 & \cdots   \\
   126 & 70 & 35 & 15 & 5 & 1 & \cdots   \\
   \vdots  & \vdots  & \vdots  & \vdots  & \vdots  & \vdots  & \ddots   \\
\end{matrix} \right)\left| {{e}^{x}} \right|.$$

We introduce ordinary numerator polynomials similar to the polynomials $^{B}{{\varphi }_{n}}\left( x \right)$. Denote
$$[n,\searrow ]\left( {{\left( xa\left( x \right) \right)}^{\prime }},xa\left( x \right) \right)=\frac{^{B}{{\alpha }_{n}}\left( x \right)}{{{\left( 1-x \right)}^{n+1}}},
\quad\left[ n,\searrow  \right]\left( 1+x{{\left( \log a\left( x \right) \right)}^{\prime }},x{{a}^{-1}}\left( x \right) \right)=\frac{^{B}\alpha _{n}^{\left( -1 \right)}\left( x \right)}{{{\left( 1-x \right)}^{n+1}}}.$$
Then
$${}^{B}\alpha _{n}^{\left( -1 \right)}\left( x \right)={{\left( -1 \right)}^{n}}{{J}_{n}}{}^{B}{{\alpha }_{n}}\left( x \right).$$
{\bfseries Example 3.} Let  ${{\tilde{\varphi }}_{n}}\left( x \right)|a\left( x \right)$, $^{B}{{\alpha }_{n}}\left( x \right)|a\left( x \right)$ denotes respectively polynomials ${{\tilde{\varphi }}_{n}}\left( x \right)$, $^{B}{{\alpha }_{n}}\left( x \right)$,  associated with the matrices ${{\left( {{\left( xa\left( x \right) \right)}^{\prime }},xa\left( x \right) \right)}_{{{e}^{x}}}}$, $\left( {{\left( xa\left( x \right) \right)}^{\prime }},xa\left( x \right) \right)$. Since
$${{\tilde{\varphi }}_{n}}\left( x \right)|C\left( x \right)=\frac{\left( 2n \right)!}{n!},$$
then
$$^{B}{{\alpha }_{n}}\left( x \right)|C\left( x \right)=\frac{\left( 2n \right)!}{n!}S_{n}^{-1}{{x}^{0}}=\sum\limits_{m=0}^{n}{\left( \begin{matrix}
   -n  \\
   m  \\
\end{matrix} \right)\left( \begin{matrix}
   2n  \\
   n-m  \\
\end{matrix} \right){{x}^{m}}};$$
$$\left( {{\left( xC\left( x \right) \right)}^{\prime }},xC\left( x \right) \right)=\left( \begin{matrix}
   1 & 0 & 0 & 0 & 0 & 0 & \cdots   \\
   2 & 1 & 0 & 0 & 0 & 0 & \cdots   \\
   6 & 3 & 1 & 0 & 0 & 0 & \cdots   \\
   20 & 10 & 4 & 1 & 0 & 0 & \cdots   \\
   70 & 35 & 15 & 5 & 1 & 0 & \cdots   \\
   252 & 126 & 56 & 21 & 6 & 1 & \cdots   \\
   \vdots  & \vdots  & \vdots  & \vdots  & \vdots  & \vdots  & \ddots   \\
\end{matrix} \right).$$
Respectively, polynomials 
$${{\left( -1 \right)}^{n}}\sum\limits_{m=0}^{n}{\left( \begin{matrix}
   2n  \\
   m  \\
\end{matrix} \right)\left( \begin{matrix}
   -n  \\
   n-m  \\
\end{matrix} \right){{x}^{m}}}$$
 are the numerator polynomials of the matrix
$$\left( 1+x{{\left( \operatorname{logC}\left( x \right) \right)}^{\prime }},x{{C}^{-1}}\left( x \right) \right)=\left( \begin{matrix}
   1 & 0 & 0 & 0 & 0 & 0 & \cdots   \\
   1 & 1 & 0 & 0 & 0 & 0 & \cdots   \\
   3 & 0 & 1 & 0 & 0 & 0 & \cdots   \\
   10 & 1 & -1 & 1 & 0 & 0 & \cdots   \\
   35 & 4 & 0 & -2 & 1 & 0 & \cdots   \\
   126 & 15 & 1 & 0 & -3 & 1 & \cdots   \\
   \vdots  & \vdots  & \vdots  & \vdots  & \vdots  & \vdots  & \ddots   \\
\end{matrix} \right).$$
{\bfseries Example 4.} Since
$${{\tilde{\varphi }}_{n}}\left( x \right)|1+x=\frac{\left( 2n \right)!}{n!}{{x}^{n-1}},$$
then
$$^{B}{{\alpha }_{n}}\left( x \right)|1+x=\frac{\left( 2n \right)!}{n!}S_{n}^{-1}{{x}^{n-1}}=\left( 2-x \right){{x}^{n-1}};$$
$$\left( {{\left( x\left( 1+x \right) \right)}^{\prime }},x\left( 1+x \right) \right)=\left( \begin{matrix}
   1 & 0 & 0 & 0 & 0 & 0 & 0 & \cdots   \\
   2 & 1 & 0 & 0 & 0 & 0 & 0 & \cdots   \\
   0 & 3 & 1 & 0 & 0 & 0 & 0 & \cdots   \\
   0 & 2 & 4 & 1 & 0 & 0 & 0 & \cdots   \\
   0 & 0 & 5 & 5 & 1 & 0 & 0 & \cdots   \\
   0 & 0 & 2 & 9 & 6 & 1 & 0 & \cdots   \\
   0 & 0 & 0 & 7 & 14 & 7 & 1 & \cdots   \\
   \vdots  & \vdots  & \vdots  & \vdots  & \vdots  & \vdots  & \vdots  & \ddots   \\
\end{matrix} \right).$$
Respectively, polynomials ${{\left( -1 \right)}^{n}}\left( -1+2x \right)$ are the numerator polynomials of the matrix
$$\left( 1+x{{\left( \log \left( 1+x \right) \right)}^{\prime }},x{{\left( 1+x \right)}^{-1}} \right)=\left( \begin{matrix}
   1 & 0 & 0 & 0 & 0 & 0 & 0 & \cdots   \\
   1 & 1 & 0 & 0 & 0 & 0 & 0 & \cdots   \\
   -1 & 0 & 1 & 0 & 0 & 0 & 0 & \cdots   \\
   1 & -1 & -1 & 1 & 0 & 0 & 0 & \cdots   \\
   -1 & 2 & 0 & -2 & 1 & 0 & 0 & \cdots   \\
   1 & -3 & 2 & 2 & -3 & 1 & 0 & \cdots   \\
   -1 & 4 & -5 & 0 & 5 & -4 & 1 & \cdots   \\
   \vdots  & \vdots  & \vdots  & \vdots  & \vdots  & \vdots  & \vdots  & \ddots   \\
\end{matrix} \right).$$
{\bfseries Example 5.} Since ${{\left( x{{\left( 1-x \right)}^{-1}} \right)}^{\prime }}={{\left( 1-x \right)}^{-2}}$, then
$$\left[ n,\searrow  \right]\left( {{\left( \frac{x}{1-x} \right)}^{\prime }},\frac{x}{1-x} \right)=\frac{1}{x}\left( \frac{1}{{{\left( 1-x \right)}^{n+1}}}-1 \right), \qquad^{B}{{\alpha }_{n}}\left( x \right)=\frac{1-{{\left( 1-x \right)}^{n+1}}}{x}.$$
Respectively, numerator polynomials of the matrix $\left( {{\left( 1-x \right)}^{-1}},x\left( 1-x \right) \right)$ are the polynomials  ${{\left( 1-x \right)}^{n+1}}+{{\left( -x \right)}^{n}}$.\\
{\bfseries Example 6.} Since ${{\left( x{{e}^{x}} \right)}^{\prime }}=\left( 1+x \right){{e}^{x}}$, then
$$\left[ n,\searrow  \right]\left( {{\left( x{{e}^{x}} \right)}^{\prime }},x{{e}^{x}} \right)=\sum\limits_{m=0}^{\infty }{\frac{{{\left( m+1 \right)}^{n}}+n{{\left( m+1 \right)}^{n-1}}}{n!}}{{x}^{m}},$$
$$^{B}{{\alpha }_{n}}\left( x \right)|{{e}^{x}}=\frac{1}{n!}\left( {{{\tilde{A}}}_{n}}\left( x \right)+n\left( 1-x \right){{{\tilde{A}}}_{n-1}}\left( x \right) \right),  \qquad{{\tilde{A}}_{0}}\left( x \right)={{A}_{0}}\left( x \right).$$
Respectively, numerator polynomials of the matrix $\left( 1+x,x{{e}^{-x}} \right)$ are the polynomials  
$$\frac{{{\left( -1 \right)}^{n}}}{n!}\left( {{A}_{n}}\left( x \right)-n\left( 1-x \right){{A}_{n-1}}\left( x \right) \right).$$
Note that
$${{\tilde{\varphi }}_{n}}\left( x \right)|{{e}^{x}}={{F}_{n}}\left( x+n+1 \right){{\left( x+1 \right)}^{n-1}}={}^{B}{{F}_{n}}\left( x+n \right){{x}^{n-1}}.$$
In general case
$${{\tilde{\varphi }}_{n}}\left( x \right)={{F}_{n}}\left( x+n+1 \right){{\tilde{u}}_{n}}\left( x+1 \right)={}^{B}{{F}_{n}}\left( x+n \right){{\tilde{u}}_{n}}\left( x \right),$$
or
$${{F}_{n}}\left( x+n+1,x \right)E{{I}_{n-1}}={}^{B}{{F}_{n}}\left( x+n,x \right){{I}_{n-1}}={{\tilde{F}}_{n}}.$$
Hence,
$$\left[ {{x}^{n}} \right]{{F}_{n}}{{x}^{p}}\left( x+n+1 \right)=\left[ {{x}^{n}} \right]{}^{B}{{F}_{n}}{{x}^{p}}\left( x+n \right)=0, \qquad p<n.$$
Here the identities for the $n$th elements of the columns of the matrices ${{F}_{n}}$, ${}^{B}{{F}_{n}}$ are manifested:
$$\sum\limits_{m=0}^{n}{{{\left( -1 \right)}^{n-m}}\left( \begin{matrix}
   2n+1  \\
   n-m  \\
\end{matrix} \right){{m}^{p}}}\left( \begin{matrix}
   m+n  \\
   n  \\
\end{matrix} \right)={{\left( -1 \right)}^{n+p}}{{\left( n+1 \right)}^{p}},$$
$$\sum\limits_{m=0}^{n}{{{\left( -1 \right)}^{n-m}}}\left( \begin{matrix}
   2n+1  \\
   n-m  \\
\end{matrix} \right){{\left( m+1 \right)}^{p}}\left( \begin{matrix}
   m+n  \\
   n  \\
\end{matrix} \right)={{\left( -1 \right)}^{n+p}}{{n}^{p}},  \qquad p\le n.$$

\section{Numerator polynomials and generalized Lagrange series}
We return to the series $_{\left( \beta  \right)}a\left( x \right)$, $_{\left( 0 \right)}a\left( x \right)=a\left( x \right)$, from Section 2. Parameter $\beta $ is defined by the identity
$$a\left( x{}_{\left( \beta  \right)}{{a}^{\beta }}\left( x \right) \right)={}_{\left( \beta  \right)}a\left( x \right),$$
so that
$${{a}^{\beta }}\left( x{}_{\left( \beta  \right)}{{a}^{\beta }}\left( x \right) \right)={}_{\left( \beta  \right)}{{a}^{\beta }}\left( x \right).$$
Denote ${{a}^{\beta }}\left( x \right)=c\left( x \right)$, ${}_{\left( \beta  \right)}{{a}^{\beta }}\left( x \right)=d\left( x \right)$. By the Lagrange inversion theorem, if
$$b\left( x \right){{c}^{\varphi }}\left( x \right)=\sum\limits_{n=0}^{\infty }{{{f}_{n}}}\left( \varphi  \right){{x}^{n}},$$
where ${{f}_{n}}\left( x \right)$ are the polynomials, then
$$b\left( xd\left( x \right) \right)\left( 1+x{{\left( \log d\left( x \right) \right)}^{\prime }} \right){{d}^{\varphi }}\left( x \right)=\sum\limits_{n=0}^{\infty }{{{f}_{n}}\left( \varphi +n \right)}{{x}^{n}}.$$
Here
$${{f}_{n}}\left( x \right)=\frac{{{s}_{n}}\left( \beta x \right)}{n!}, \qquad{{s}_{n}}\left( x \right)=\left[ n,\to  \right]{{\left( b\left( x \right),\log a\left( x \right) \right)}_{{{e}^{x}}}}.$$
Thus,
$${{s}_{n}}\left( \beta x+\beta n \right)=\left[ n,\to  \right]{{\left( b\left( x{}_{\left( \beta  \right)}{{a}^{\beta }}\left( x \right) \right)\left( 1+x{{\left( \log {}_{\left( \beta  \right)}{{a}^{\beta }}\left( x \right) \right)}^{\prime }} \right),\log {}_{\left( \beta  \right)}{{a}^{\beta }}\left( x \right) \right)}_{{{e}^{x}}}},$$ 
$${{s}_{n}}\left( x+\beta n \right)=\left[ n,\to  \right]{{\left( b\left( x{}_{\left( \beta  \right)}{{a}^{\beta }}\left( x \right) \right)\left( 1+x{{\left( \log {}_{\left( \beta  \right)}{{a}^{\beta }}\left( x \right) \right)}^{\prime }} \right),\log {}_{\left( \beta  \right)}a\left( x \right) \right)}_{{{e}^{x}}}}.$$
Denote
$$\left[ n,\searrow  \right]\left( b\left( x{}_{\left( \beta  \right)}{{a}^{\beta }}\left( x \right) \right)\left( 1+x\beta {{\left( \log {}_{\left( \beta  \right)}a\left( x \right) \right)}^{\prime }} \right),x{}_{\left( \beta  \right)}a\left( x \right) \right)=\frac{_{\left( \beta  \right)}{{g}_{n}}\left( x \right)}{{{\left( 1-x \right)}^{n+1}}}.$$
We introduce the matrices $G_{n}^{\beta }={{U}_{n}}{{E}^{n\beta }}U_{n}^{-1}$. For example,
$${{G}_{2}}=\left( \begin{matrix}
   6 & 3 & 1  \\
   -8 & -3 & 0  \\
   3 & 1 & 0  \\
\end{matrix} \right),   \ {{G}_{3}}=\left( \begin{matrix}
   20 & 10 & 4 & 1  \\
   -45 & -20 & -6 & 0  \\
   36 & 15 & 4 & 0  \\
   -10 & -4 & -1 & 0  \\
\end{matrix} \right),\
{{G}_{4}}=\left( \begin{matrix}
   70 & 35 & 15 & 5 & 1  \\
   -224 & -105 & -40 & -10 & 0  \\
   280 & 126 & 45 & 10 & 0  \\
   -160 & -70 & -24 & -5 & 0  \\
   35 & 15 & 5 & 1 & 0  \\
\end{matrix} \right).$$
Then
$$G_{n}^{\beta }{{g}_{n}}\left( x \right)={}_{\left( \beta  \right)}{{g}_{n}}\left( x \right).$$
{\bfseries Theorem 15.}
$$G_{n}^{-\beta }={{J}_{n}}G_{n}^{\beta }{{J}_{n}}.$$
{\bfseries Proof.} Since $E\left( 1,-x \right)=\left( 1,-x \right){{E}^{-1}}$, by Theorem 6
$${{J}_{n}}{{U}_{n}}{{E}^{n\beta }}U_{n}^{-1}{{J}_{n}}={{U}_{n}}E\left( 1,-x \right){{E}^{n\beta }}E\left( 1,-x \right)U_{n}^{-1}=$$
$$={{U}_{n}}\left( 1,-x \right){{E}^{n\beta }}\left( 1,-x \right)U_{n}^{-1}={{U}_{n}}{{E}^{-n\beta }}U_{n}^{-1}.$$
Thus,
$$G_{2}^{-1}=\left( \begin{matrix}
   0 & 1 & 3  \\
   0 & -3 & -8  \\
   1 & 3 & 6  \\
\end{matrix} \right),\  G_{3}^{-1}=\left( \begin{matrix}
   0 & -1 & -4 & -10  \\
   0 & 4 & 15 & 36  \\
   0 & -6 & -20 & -45  \\
   1 & 4 & 10 & 20  \\
\end{matrix} \right),
\ G_{4}^{-1}=\left( \begin{matrix}
   0 & 1 & 5 & 15 & 35  \\
   0 & -5 & -24 & -70 & -160  \\
   0 & 10 & 45 & 126 & 280  \\
   0 & -10 & -40 & -105 & -224  \\
   1 & 5 & 15 & 35 & 70  \\
\end{matrix} \right).$$
{\bfseries Theorem 16.}
$$G_{n}^{\beta }=V_{n}^{-1}{{\left( {{\left( 1+x \right)}^{n\beta }},x \right)}^{T}}{{V}_{n}}.$$
{\bfseries Proof.} Since
$$n!\left| {{e}^{x}} \right|{{V}_{n}}{{U}_{n}}{{x}^{p}}=\left[ p,\to  \right]{{\left( 1,{{e}^{x}}-1 \right)}_{{{e}^{x}}}},$$
$$\left( {1}/{n!}\; \right)U_{n}^{-1}V_{n}^{-1}{{\left| {{e}^{x}} \right|}^{-1}}{{x}^{p}}=\left[ p,\to  \right]{{\left( 1,\log \left( 1+x \right) \right)}_{{{e}^{x}}}},$$
$${{\left( 1,\log \left( 1+x \right) \right)}_{{{e}^{x}}}}{{\left( {{e}^{n\beta }},x \right)}_{{{e}^{x}}}}{{\left( 1,{{e}^{x}}-1 \right)}_{{{e}^{x}}}}={{\left( {{\left( 1+x \right)}^{n\beta }},x \right)}_{{{e}^{x}}}},$$
then
$${{V}_{n}}{{U}_{n}}{{E}^{n\beta }}U_{n}^{-1}V_{n}^{-1}={{\left( {{\left( 1+x \right)}^{n\beta }},x \right)}^{T}}{{I}_{n}}.$$

For example,
$${{G}_{3}}=\left( \begin{matrix}
   1 & 0 & 0 & 0  \\
   -3 & 1 & 0 & 0  \\
   3 & -2 & 1 & 0  \\
   -1 & 1 & -1 & 1  \\
\end{matrix} \right)\left( \begin{matrix}
   1 & 3 & 3 & 1  \\
   0 & 1 & 3 & 3  \\
   0 & 0 & 1 & 3  \\
   0 & 0 & 0 & 1  \\
\end{matrix} \right)\left( \begin{matrix}
   1 & 0 & 0 & 0  \\
   3 & 1 & 0 & 0  \\
   3 & 2 & 1 & 0  \\
   1 & 1 & 1 & 1  \\
\end{matrix} \right).$$
{\bfseries Theorem 17.}
$$G_{n}^{\beta }{{x}^{p}}=\sum\limits_{m=0}^{n}{\left( \begin{matrix}
   -n\beta +p  \\
   m  \\
\end{matrix} \right)\left( \begin{matrix}
   n\beta +n-p  \\
   n-m  \\
\end{matrix} \right){{x}^{m}}}.$$
{\bfseries Proof.} Taking into account the polynomial argument, we prove the theorem for positive integers $\beta $. 
$${{\left( {{\left( 1+x \right)}^{n\beta }},x \right)}^{T}}{{V}_{n}}{{x}^{p}}=\sum\limits_{i=0}^{n}{\left( \begin{matrix}
   n\beta +n-p  \\
   n-i  \\
\end{matrix} \right){{x}^{i}}},$$
$$\left[ {{x}^{m}} \right]V_{n}^{-1}{{\left( {{\left( 1+x \right)}^{n\beta }},x \right)}^{T}}{{V}_{n}}{{x}^{p}}=$$
$$=\sum\limits_{i=0}^{m}{{{\left( -1 \right)}^{m-i}}}\left( \begin{matrix}
   n-i  \\
   m-i  \\
\end{matrix} \right)\left( \begin{matrix}
   n\beta +n-p  \\
   n-i  \\
\end{matrix} \right)\frac{\left( n\beta +m-p \right)!}{\left( n\beta +m-p \right)!}=$$ 
$$=\left( \begin{matrix}
   n\beta +n-p  \\
   n-m  \\
\end{matrix} \right)\sum\limits_{i=0}^{m}{{{\left( -1 \right)}^{m-i}}}\left( \begin{matrix}
   n\beta +m-p  \\
   m-i  \\
\end{matrix} \right)=$$
$$=\left( \begin{matrix}
   n\beta +n-p  \\
   n-m  \\
\end{matrix} \right){{\left( -1 \right)}^{m}}\left( \begin{matrix}
   n\beta +m-p-1  \\
   m  \\
\end{matrix} \right)=\left( \begin{matrix}
   n\beta +n-p  \\
   n-m  \\
\end{matrix} \right)\left( \begin{matrix}
   -n\beta +p  \\
   m  \\
\end{matrix} \right).$$

We introduce the matrices ${{X}_{n}}=V_{n}^{-1}{{\left( x,x \right)}^{T}}{{V}_{n}}$. Since $V_{n}^{-1}=\left( {{\left( 1-x \right)}^{n+1}},x \right)P{{I}_{n}}$,  we find:
$${{X}_{n}}{{x}^{0}}=\frac{1-x-{{\left( 1-x \right)}^{n+1}}}{x} ,   \qquad{{X}_{n}}{{x}^{p}}={{x}^{p-1}}\left( 1-x \right).$$
Then
$$G_{n}^{\beta }={{\left( {{I}_{n}}+{{X}_{n}} \right)}^{n\beta }}=\sum\limits_{m=0}^{n}{\left( \begin{matrix}
   n\beta   \\
   m  \\
\end{matrix} \right)X_{n}^{m}}.$$
For example,
$${{G}_{3}}={{I}_{3}}+3\left( \begin{matrix}
   3 & 1 & 0 & 0  \\
   -6 & -1 & 1 & 0  \\
   4 & 0 & -1 & 1  \\
   -1 & 0 & 0 & -1  \\
\end{matrix} \right)+3\left( \begin{matrix}
   3 & 2 & 1 & 0  \\
   -8 & -5 & -2 & 1  \\
   7 & 4 & 1 & -2  \\
   -2 & -1 & 0 & 1  \\
\end{matrix} \right)+\left( \begin{matrix}
   1 & 1 & 1 & 1  \\
   -3 & -3 & -3 & -3  \\
   3 & 3 & 3 & 3  \\
   -1 & -1 & -1 & -1  \\
\end{matrix} \right),$$
$$G_{3}^{-1}={{I}_{3}}-3\left( \begin{matrix}
   3 & 1 & 0 & 0  \\
   -6 & -1 & 1 & 0  \\
   4 & 0 & -1 & 1  \\
   -1 & 0 & 0 & -1  \\
\end{matrix} \right)+6\left( \begin{matrix}
   3 & 2 & 1 & 0  \\
   -8 & -5 & -2 & 1  \\
   7 & 4 & 1 & -2  \\
   -2 & -1 & 0 & 1  \\
\end{matrix} \right)-10\left( \begin{matrix}
   1 & 1 & 1 & 1  \\
   -3 & -3 & -3 & -3  \\
   3 & 3 & 3 & 3  \\
   -1 & -1 & -1 & -1  \\
\end{matrix} \right).$$
Thus,
$${{I}_{n}}+{{X}_{n}}=G_{n}^{{1}/{n}\;}={{U}_{n}}EU_{n}^{-1}.$$
For example,
$$G_{2}^{{1}/{2}\;}=\left( \begin{matrix}
   3 & 1 & 0  \\
   -3 & 0 & 1  \\
   1 & 0 & 0  \\
\end{matrix} \right),   \qquad G_{3}^{{1}/{3}\;}\left( \begin{matrix}
   4 & 1 & 0 & 0  \\
   -6 & 0 & 1 & 0  \\
   4 & 0 & 0 & 1  \\
   -1 & 0 & 0 & 0  \\
\end{matrix} \right),  \qquad G_{4}^{{1}/{4}\;}=\left( \begin{matrix}
   5 & 1 & 0 & 0 & 0  \\
   -10 & 0 & 1 & 0 & 0  \\
   10 & 0 & 0 & 1 & 0  \\
   -5 & 0 & 0 & 0 & 1  \\
   1 & 0 & 0 & 0 & 0  \\
\end{matrix} \right),$$
$$G_{2}^{{-1}/{2}\;}=\left( \begin{matrix}
   0 & 0 & 1  \\
   1 & 0 & -3  \\
   0 & 1 & 3  \\
\end{matrix} \right),  \quad G_{3}^{{-1}/{3}\;}\left( \begin{matrix}
   0 & 0 & 0 & -1  \\
   1 & 0 & 0 & 4  \\
   0 & 1 & 0 & -6  \\
   0 & 0 & 1 & 4  \\
\end{matrix} \right),  \quad G_{4}^{{-1}/{4}\;}=\left( \begin{matrix}
   0 & 0 & 0 & 0 & 1  \\
   1 & 0 & 0 & 0 & -5  \\
   0 & 1 & 0 & 0 & 10  \\
   0 & 0 & 1 & 0 & -10  \\
   0 & 0 & 0 & 1 & 5  \\
\end{matrix} \right).$$
From Theorem 10 it follows that
$$\left( {{\left( 1-x \right)}^{-m}},x \right)G_{n}^{\beta }\left( {{\left( 1-x \right)}^{m}},x \right){{I}_{n-m}}=G_{n-m}^{\frac{n\beta }{n-m}}.$$
In particular,
$$\left( {{\left( 1-x \right)}^{-m}},x \right)G_{n}^{{1}/{n}\;}\left( {{\left( 1-x \right)}^{m}},x \right){{I}_{n-m}}=G_{n-m}^{{1}/{\left( n-m \right)}\;}.$$
{\bfseries Example 7.} Let $b\left( x \right)=1$, $a\left( x \right)=1+x$, $_{\left( \beta  \right)}a\left( x \right)$ is the generalized binomial series. Then polynomials $_{\left( \beta  \right)}{{g}_{n}}\left( x \right)=G_{n}^{\beta }{{x}^{n}}$ are the numerator polynomials of the matrix
$$\left( 1+x{{\left( \log {}_{\left( \beta  \right)}{{a}^{\beta }}\left( x \right) \right)}^{\prime }},{{x}_{\left( \beta  \right)}}a\left( x \right) \right).$$
But since $_{\left( 1 \right)}{{g}_{n}}\left( x \right)=1$, $G_{n}^{\beta }{{x}^{0}}=G_{n}^{\beta +1}{{x}^{n}}$, then this matrix can also be represented in the form
$$\left( _{\left( \beta  \right)}a\left( x \right)\left( 1+x{{\left( \log {}_{\left( \beta  \right)}{{a}^{\beta -1}}\left( x \right) \right)}^{\prime }} \right),{{x}_{\left( \beta  \right)}}a\left( x \right) \right).$$
Hence, here the property of the generalized binomial series is manifested:
$$_{\left( \beta  \right)}a\left( x \right)\left( 1+x{{\left( \log {}_{\left( \beta  \right)}{{a}^{\beta -1}}\left( x \right) \right)}^{\prime }} \right)=1+x{{\left( \log {}_{\left( \beta  \right)}{{a}^{\beta }}\left( x \right) \right)}^{\prime }}.$$
{\bfseries Example 8.} Let $_{\left( \beta  \right)}a\left( x \right)$ is the generalized binomial series. Then polynomials $G_{n}^{\beta }x$ are the numerator polynomials of the matrix
$$\left( 1+x{{\left( \log {}_{\left( \beta +1 \right)}{{a}^{\beta }}\left( x \right) \right)}^{\prime }},{{x}_{\left( \beta +1 \right)}}a\left( x \right) \right),$$
polynomials $G_{n}^{\beta }{{x}^{n-1}}$ are the numerator polynomials of the matrix
$$\left( _{\left( \beta  \right)}a\left( x \right)\left( 1+x{{\left( \log {}_{\left( \beta  \right)}{{a}^{\beta }}\left( x \right) \right)}^{\prime }} \right),{{x}_{\left( \beta  \right)}}a\left( x \right) \right).$$
Since $G_{n}^{-\beta }x={{J}_{n}}G_{n}^{\beta }{{x}^{n-1}}$, then  matrix
$$\left( 1+x{{\left( \log {}_{\left( 1-\beta  \right)}{{a}^{-\beta }}\left( x \right) \right)}^{\prime }},{{x}_{\left( 1-\beta  \right)}}a\left( x \right) \right)$$
coincides with the matrix
$$\left( 1,-x \right)\left( 1+x{{\left( \log {}_{\left( \beta  \right)}{{a}^{\beta }}\left( x \right) \right)}^{\prime }},{{x}_{\left( \beta  \right)}}{{a}^{-1}}\left( x \right) \right)\left( 1,-x \right),$$
matrix
$$\left( _{\left( -\beta  \right)}a\left( x \right)\left( 1+x{{\left( \log {}_{\left( -\beta  \right)}{{a}^{-\beta }}\left( x \right) \right)}^{\prime }} \right),{{x}_{\left( -\beta  \right)}}a\left( x \right) \right)$$
coincides with the matrix
$$\left( 1,-x \right)\left( _{\left( \beta +1 \right)}{{a}^{-1}}\left( x \right)\left( 1+x{{\left( \log {}_{\left( \beta +1 \right)}{{a}^{\beta }}\left( x \right) \right)}^{\prime }} \right),{{x}_{\left( \beta +1 \right)}}{{a}^{-1}}\left( x \right) \right)\left( 1,-x \right).$$

Denote
$$\left[ n,\searrow  \right]{{\left( b\left( x{}_{\left( \beta  \right)}{{a}^{\beta }}\left( x \right) \right)\left( 1+x\beta {{\left( \log {}_{\left( \beta  \right)}a\left( x \right) \right)}^{\prime }} \right),x{}_{\left( \beta  \right)}a\left( x \right) \right)}_{{{e}^{x}}}}=\frac{_{\left( \beta  \right)}{{h}_{n}}\left( x \right)}{{{\left( 1-x \right)}^{2n+1}}}.$$
We introduce the matrices  $H_{n}^{\beta }={{F}_{n}}{{E}^{n\beta }}F_{n}^{-1}$. For example,
$${{H}_{2}}=\frac{1}{6}\left( \begin{matrix}
   15 & 5 & 1  \\
   -12 & 2 & 4  \\
   3 & -1 & 1  \\
\end{matrix} \right),  \qquad {{H}_{3}}=\frac{1}{20}\left( \begin{matrix}
   84 & 28 & 7 & 1  \\
   -108 & -4 & 15 & 9  \\
   54 & -6 & -1 & 9  \\
   -10 & 2 & -1 & 1  \\
\end{matrix} \right),$$
$$ {{H}_{4}}=\frac{1}{70}\left( \begin{matrix}
   495 & 165 & {135}/{3}\; & 9 & 1  \\
   -880 & -110 & {160}/{3}\; & 44 & 16  \\
   660 & 0 & {-90}/{3}\; & 24 & 36  \\
   -240 & 20 & 0 & -6 & 16  \\
   35 & -5 & {5}/{3}\; & -1 & 1  \\
\end{matrix} \right).$$ 
Then
$$H_{n}^{\beta }{{h}_{n}}\left( x \right)={}_{\left( \beta  \right)}{{h}_{n}}\left( x \right).$$
{\bfseries Theorem 18.}
$$H_{n}^{-\beta }={{J}_{n}}H_{n}^{\beta }{{J}_{n}}.$$
{\bfseries Proof.} By Theorem 11
$${{J}_{n}}{{F}_{n}}{{E}^{n\beta }}F_{n}^{-1}{{J}_{n}}={{F}_{n}}{{E}^{n+1}}\left( 1,-x \right){{E}^{n\beta }}{{E}^{n+1}}\left( 1,-x \right)F_{n}^{-1}=$$
$$={{F}_{n}}\left( 1,-x \right){{E}^{n\beta }}\left( 1,-x \right)F_{n}^{-1}={{F}_{n}}{{E}^{-n\beta }}F_{n}^{-1}.$$

Matrix $H_{n}^{\beta }$ can be represented in the form
$$H_{n}^{\beta }={{S}_{n}}G_{n}^{\beta }S_{n}^{-1}=V_{n}^{-1}{{C}_{n}}{{\left( {{\left( 1+x \right)}^{n\beta }},x \right)}^{T}}C_{n}^{-1}{{V}_{n}}.$$
For example,
$${{H}_{3}}=\left( \begin{matrix}
   1 & 0 & 0 & 0  \\
   -3 & 1 & 0 & 0  \\
   3 & -2 & 1 & 0  \\
   -1 & 1 & -1 & 1  \\
\end{matrix} \right)\left( \begin{matrix}
   1 & 0 & 0 & 0  \\
   0 & 4 & 0 & 0  \\
   0 & 0 & 10 & 0  \\
   0 & 0 & 0 & 20  \\
\end{matrix} \right)\left( \begin{matrix}
   1 & 3 & 3 & 1  \\
   0 & 1 & 3 & 3  \\
   0 & 0 & 1 & 3  \\
   0 & 0 & 0 & 1  \\
\end{matrix} \right)\left( \begin{matrix}
   1 & 0 & 0 & 0  \\
   0 & \frac{1}{4} & 0 & 0  \\
   0 & 0 & \frac{1}{10} & 0  \\
   0 & 0 & 0 & \frac{1}{20}  \\
\end{matrix} \right)\left( \begin{matrix}
   1 & 0 & 0 & 0  \\
   3 & 1 & 0 & 0  \\
   3 & 2 & 1 & 0  \\
   1 & 1 & 1 & 1  \\
\end{matrix} \right).$$
Denote
$${{t}_{n}}\left( \varphi |\beta ,x \right)=\sum\limits_{m=0}^{n}{\left( \begin{matrix}
   \varphi   \\
   m  \\
\end{matrix} \right)\left( \begin{matrix}
   \beta   \\
   n-m  \\
\end{matrix} \right){{x}^{m}}}.$$
{\bfseries Theorem 19.}
$$H_{n}^{\beta }{{x}^{p}}=\sum\limits_{m=p}^{n}{\left( \begin{matrix}
   n-p  \\
   n-m  \\
\end{matrix} \right){{\left( \begin{matrix}
   n+m  \\
   m  \\
\end{matrix} \right)}^{-1}}}{{\left( 1-x \right)}^{n-m}}{{t}_{m}}\left( -n\beta +n+m|n\beta ,x \right).$$
{\bfseries Proof.} Since
$$\frac{1}{n!}\left[ {{x}^{m}} \right]V_{p}^{-1}{{C}_{n}}{{\left( {{\left( 1+x \right)}^{n\beta }},x \right)}^{T}}{{x}^{p}}=$$
$$=\sum\limits_{i=0}^{m}{{{\left( -1 \right)}^{m-i}}}\left( \begin{matrix}
   p-i  \\
   m-i  \\
\end{matrix} \right)\left( \begin{matrix}
   n\beta   \\
   p-i  \\
\end{matrix} \right)\left( \begin{matrix}
   n+i  \\
   i  \\
\end{matrix} \right)\frac{\left( n\beta +m-p \right)!}{\left( n\beta +m-p \right)!}=,$$
$$=\left( \begin{matrix}
   n\beta   \\
   p-m  \\
\end{matrix} \right)\sum\limits_{i=0}^{m}{{{\left( -1 \right)}^{m-i}}\left( \begin{matrix}
   n+i  \\
   i  \\
\end{matrix} \right)}\left( \begin{matrix}
   n\beta +m-p  \\
   m-i  \\
\end{matrix} \right)=$$
$$=\left( \begin{matrix}
   n\beta   \\
   p-m  \\
\end{matrix} \right){{\left( -1 \right)}^{m}}\left( \begin{matrix}
   n\beta +m-p-n-1  \\
   m  \\
\end{matrix} \right)=\left( \begin{matrix}
   n\beta   \\
   p-m  \\
\end{matrix} \right)\left( \begin{matrix}
   -n\beta +n+p  \\
   m  \\
\end{matrix} \right),$$
then
$$V_{n}^{-1}{{C}_{n}}{{\left( {{\left( 1+x \right)}^{n\beta }},x \right)}^{T}}{{x}^{p}}=\left( {{\left( 1-x \right)}^{n-p}},x \right)V_{p}^{-1}{{C}_{n}}{{\left( {{\left( 1+x \right)}^{n\beta }},x \right)}^{T}}{{x}^{p}}=$$
$$=n!{{\left( 1-x \right)}^{n-p}}{{t}_{p}}\left( -n\beta +n+p|n\beta ,x \right).$$
It remains to add that
$$C_{n}^{-1}{{V}_{n}}{{x}^{p}}=\frac{1}{n!}\sum\limits_{m=p}^{n}{\left( \begin{matrix}
   n-p  \\
   n-m  \\
\end{matrix} \right){{\left( \begin{matrix}
   n+m  \\
   m  \\
\end{matrix} \right)}^{-1}}}{{x}^{m}}.$$

In particular,
$$H_{n}^{\beta }{{x}^{n}}={{\left( \begin{matrix}
   2n  \\
   n  \\
\end{matrix} \right)}^{-1}}\sum\limits_{m=0}^{n}{\left( \begin{matrix}
   -n\beta +2n  \\
   m  \\
\end{matrix} \right)\left( \begin{matrix}
   n\beta   \\
   n-m  \\
\end{matrix} \right){{x}^{m}}}.$$
Respectively, by Theorem 18
$$H_{n}^{\beta }{{x}^{0}}={{\left( \begin{matrix}
   2n  \\
   n  \\
\end{matrix} \right)}^{-1}}\sum\limits_{m=0}^{n}{\left( \begin{matrix}
   -n\beta   \\
   m  \\
\end{matrix} \right)\left( \begin{matrix}
   n\beta +2n  \\
   n-m  \\
\end{matrix} \right){{x}^{m}}}.$$
{\bfseries Example 9.} Let $_{\left( \beta  \right)}a\left( x \right)$ is the generalized binomial series. Then polynomials $_{\left( \beta  \right)}{{h}_{n}}\left( x \right)=\frac{\left( 2n \right)!}{n!}H_{n}^{\beta }{{x}^{n}}$  are the numerator polynomials of the matrix
$${{\left( 1+x{{\left( \log {}_{\left( \beta  \right)}{{a}^{\beta }}\left( x \right) \right)}^{\prime }},{{x}_{\left( \beta  \right)}}a\left( x \right) \right)}_{{{e}^{x}}}}.$$
Since $_{\left( 2 \right)}{{h}_{n}}\left( x \right)=\frac{\left( 2n \right)!}{n!}$, then polynomials $\frac{\left( 2n \right)!}{n!}H_{n}^{\beta }{{x}^{0}}=\frac{\left( 2n \right)!}{n!}H_{n}^{\beta +2}{{x}^{n}}$ are the numerator polynomials of the matrix
$${{\left( 1+x{{\left( \log {}_{\left( \beta +2 \right)}{{a}^{\beta +2}}\left( x \right) \right)}^{\prime }},{{x}_{\left( \beta +2 \right)}}a\left( x \right) \right)}_{{{e}^{x}}}}.$$
Since 
$$H_{n}^{-\beta }{{x}^{0}}={{J}_{n}}H_{n}^{\beta }{{x}^{n}},   \qquad{{\left( 1,x{}_{\left( \beta  \right)}a\left( x \right) \right)}^{-1}}=\left( 1,{{x}_{\left( \beta -1 \right)}}{{a}^{-1}}\left( x \right) \right),$$
$${{\left( {{x}_{\left( \beta -1 \right)}}{{a}^{-1}}\left( x \right) \right)}^{\prime }}\left( 1,{{x}_{\left( \beta -1 \right)}}{{a}^{-1}}\left( x \right) \right)\left( 1+x{{\left( \log {}_{\left( \beta  \right)}{{a}^{\beta }}\left( x \right) \right)}^{\prime }} \right)=$$
$$=\left( 1+x{{\left( \log {}_{\left( \beta -1 \right)}{{a}^{\beta -1}}\left( x \right) \right)}^{\prime }} \right){}_{\left( \beta -1 \right)}{{a}^{-1}}\left( x \right),$$
 then matrix
$$\left( 1+x{{\left( \log {}_{\left( 2-\beta  \right)}{{a}^{2-\beta }}\left( x \right) \right)}^{\prime }},{{x}_{\left( 2-\beta  \right)}}a\left( x \right) \right)$$
coincides with the matrix
$$\left( 1,-x \right)\left( \left( 1+x{{\left( \log {}_{\left( \beta -1 \right)}{{a}^{\beta -1}}\left( x \right) \right)}^{\prime }} \right){}_{\left( \beta -1 \right)}{{a}^{-1}}\left( x \right),{{x}_{\left( \beta -1 \right)}}{{a}^{-1}}\left( x \right) \right)\left( 1,-x \right).$$

Denote
$$\left[ n,\searrow  \right]{{\left( 1,x{}_{\left( \beta  \right)}a\left( x \right) \right)}_{{{e}^{x}}}}=\frac{_{\left( \beta  \right)}{{\varphi }_{n}}\left( x \right)}{{{\left( 1-x \right)}^{2n+1}}}, \qquad\frac{1}{x}{}_{\left( \beta  \right)}{{\varphi }_{n}}\left( x \right)={}_{\left( \beta  \right)}{{\tilde{\varphi }}_{n}}\left( x \right).$$
We introduce the matrices $T_{n}^{\beta }={{\tilde{F}}_{n}}{{E}^{n\beta }}\tilde{F}_{n}^{-1}$. For example,
$${{T}_{2}}=\frac{1}{2}\left( \begin{matrix}
   3 & 1  \\
   -1 & 1  \\
\end{matrix} \right),  \qquad{{T}_{3}}=\frac{1}{5}\left( \begin{matrix}
   12 & 4 & 1  \\
   -9 & 2 & 3  \\
   2 & -1 & 1  \\
\end{matrix} \right),  \qquad{{T}_{4}}=\frac{1}{14}\left( \begin{matrix}
   55 & {55}/{3}\; & 5 & 1  \\
   -66 & 0 & 10 & 6  \\
   30 & -6 & 0 & 6  \\
   -5 & {5}/{3}\; & -1 & 1  \\
\end{matrix} \right).$$
Then 
$$T_{n}^{\beta }{{\tilde{\varphi }}_{n}}\left( x \right)={}_{\left( \beta  \right)}{{\tilde{\varphi }}_{n}}\left( x \right).$$
{\bfseries Theorem 20.}
$$T_{n}^{-\beta }={{\tilde{J}}_{n}}T_{n}^{\beta }{{\tilde{J}}_{n}}.$$
{\bfseries Proof.} By Theorem 6
$${{\tilde{J}}_{n}}{{\tilde{F}}_{n}}{{E}^{n\beta }}\tilde{F}_{n}^{-1}{{\tilde{J}}_{n}}={{\tilde{F}}_{n}}{{E}^{n}}\left( 1,-x \right){{E}^{n\beta }}{{E}^{n}}\left( 1,-x \right)\tilde{F}_{n}^{-1}=$$
$$={{\tilde{F}}_{n}}\left( 1,-x \right){{E}^{n\beta }}\left( 1,-x \right)\tilde{F}_{n}^{-1}={{\tilde{F}}_{n}}{{E}^{-n\beta }}\tilde{F}_{n}^{-1}.$$

Matrix  $T_{n}^{\beta }$ can be represented in the form
$$T_{n}^{\beta }={{\tilde{S}}_{n}}A_{n}^{\beta }\tilde{S}_{n}^{-1}=\tilde{V}_{n}^{-1}{{\tilde{C}}_{n}}\tilde{D}{{\left( {{\left( 1+x \right)}^{n\beta }},x \right)}^{T}}{{\tilde{D}}^{-1}}\tilde{C}_{n}^{-1}{{\tilde{V}}_{n}},$$  
where
$${{\tilde{C}}_{n}}\tilde{D}{{x}^{p}}=\left( n+1 \right)!\left( \begin{matrix}
   n+1+p  \\
   p  \\
\end{matrix} \right){{x}^{p}}.$$
For example,
$${{T}_{4}}=\left( \begin{matrix}
   1 & 0 & 0 & 0  \\
   -3 & 1 & 0 & 0  \\
   3 & -2 & 1 & 0  \\
   -1 & 1 & -1 & 1  \\
\end{matrix} \right)\left( \begin{matrix}
   1 & 0 & 0 & 0  \\
   0 & 6 & 0 & 0  \\
   0 & 0 & 21 & 0  \\
   0 & 0 & 0 & 56  \\
\end{matrix} \right)\left( \begin{matrix}
   1 & 4 & 6 & 4  \\
   0 & 1 & 4 & 6  \\
   0 & 0 & 1 & 4  \\
   0 & 0 & 0 & 1  \\
\end{matrix} \right)\left( \begin{matrix}
   1 & 0 & 0 & 0  \\
   0 & \frac{1}{6} & 0 & 0  \\
   0 & 0 & \frac{1}{21} & 0  \\
   0 & 0 & 0 & \frac{1}{56}  \\
\end{matrix} \right)\left( \begin{matrix}
   1 & 0 & 0 & 0  \\
   3 & 1 & 0 & 0  \\
   3 & 2 & 1 & 0  \\
   1 & 1 & 1 & 1  \\
\end{matrix} \right).$$
{\bfseries Theorem 21.}
$$T_{n}^{\beta }{{x}^{p}}=\sum\limits_{m=p}^{n-1}{\left( \begin{matrix}
   n-1-p  \\
   n-1-m  \\
\end{matrix} \right){{\left( \begin{matrix}
   n+1+m  \\
   m  \\
\end{matrix} \right)}^{-1}}}{{\left( 1-x \right)}^{n-m-1}}{{t}_{m}}\left( -n\beta +n+m+1|n\beta ,x \right).$$
{\bfseries Proof.} In this case $p=0$, $1$, … , $n-1$. Since
$$\frac{1}{\left( n+1 \right)!}\left[ {{x}^{m}} \right]\tilde{V}_{p+1}^{-1}{{\tilde{C}}_{n}}\tilde{D}{{\left( {{\left( 1+x \right)}^{n\beta }},x \right)}^{T}}{{x}^{p}}=$$
$$=\sum\limits_{i=0}^{m}{{{\left( -1 \right)}^{m-i}}}\left( \begin{matrix}
   p-i  \\
   m-i  \\
\end{matrix} \right)\left( \begin{matrix}
   n\beta   \\
   p-i  \\
\end{matrix} \right)\left( \begin{matrix}
   n+1+i  \\
   i  \\
\end{matrix} \right)\frac{\left( n\beta +m-p \right)!}{\left( n\beta +m-p \right)!}=$$
$$=\left( \begin{matrix}
   n\beta   \\
   p-m  \\
\end{matrix} \right)\sum\limits_{i=0}^{m}{{{\left( -1 \right)}^{m-i}}\left( \begin{matrix}
   n+1+i  \\
   i  \\
\end{matrix} \right)}\left( \begin{matrix}
   n\beta +m-p  \\
   m-i  \\
\end{matrix} \right)=$$
$$=\left( \begin{matrix}
   n\beta   \\
   p-m  \\
\end{matrix} \right){{\left( -1 \right)}^{m}}\left( \begin{matrix}
   n\beta +m-p-n-2  \\
   m  \\
\end{matrix} \right)=\left( \begin{matrix}
   n\beta   \\
   p-m  \\
\end{matrix} \right)\left( \begin{matrix}
   -n\beta +n+p+1  \\
   m  \\
\end{matrix} \right),$$
then
$$\tilde{V}_{n}^{-1}{{\tilde{C}}_{n}}\tilde{D}{{\left( {{\left( 1+x \right)}^{n\beta }},x \right)}^{T}}{{x}^{p}}=\left( {{\left( 1-x \right)}^{n-p-1}},x \right)\tilde{V}_{p+1}^{-1}{{\tilde{C}}_{n}}\tilde{D}{{\left( {{\left( 1+x \right)}^{n\beta }},x \right)}^{T}}{{x}^{p}}=$$
$$=\left( n+1 \right)!{{\left( 1-x \right)}^{n-p-1}}{{t}_{p}}\left( -n\beta +n+p+1|n\beta ,x \right).$$
It remains to add that
$${{\tilde{D}}^{-1}}\tilde{C}_{n}^{-1}{{\tilde{V}}_{n}}{{x}^{p}}=\frac{1}{\left( n+1 \right)!}\sum\limits_{m=p}^{n-1}{\left( \begin{matrix}
   n-1-p  \\
   n-1-m  \\
\end{matrix} \right){{\left( \begin{matrix}
   n+1+m  \\
   m  \\
\end{matrix} \right)}^{-1}}}{{x}^{m}}.$$

In particular,
$$T_{n}^{\beta }{{x}^{n-1}}={{\left( \begin{matrix}
   2n  \\
   n-1  \\
\end{matrix} \right)}^{-1}}\sum\limits_{m=0}^{n-1}{\left( \begin{matrix}
   n\left( 2-\beta  \right)  \\
   m  \\
\end{matrix} \right)\left( \begin{matrix}
   n\beta   \\
   n-1-m  \\
\end{matrix} \right){{x}^{m}}}.$$
Respectively, by Theorem 20
$$T_{n}^{\beta }{{x}^{0}}={{\left( \begin{matrix}
   2n  \\
   n-1  \\
\end{matrix} \right)}^{-1}}\sum\limits_{m=0}^{n-1}{\left( \begin{matrix}
   -n\beta   \\
   m  \\
\end{matrix} \right)\left( \begin{matrix}
   n\left( 2+\beta  \right)  \\
   n-1-m  \\
\end{matrix} \right){{x}^{m}}}.$$

E-mail: {evgeniy\symbol{"5F}burlachenko@list.ru}
\end{document}